\RequirePackage{float}
\documentclass[pdflatex,iicol,sn-mathphys-ay]{sn-jnl}
\usepackage{graphicx}
\usepackage{xcolor}
\usepackage{amsmath,amssymb,amsfonts,amsthm}
\usepackage{booktabs,array,multirow}
\usepackage{siunitx}
\usepackage{placeins}
\usepackage{algorithm}
\usepackage{diagbox}
\usepackage{adjustbox}

\usepackage{tikz}

\usetikzlibrary{svg.path}

\definecolor{orcidlogocol}{HTML}{A6CE39}
\tikzset{
  orcidlogo/.pic={
    \fill[orcidlogocol] svg{M256,128c0,70.7-57.3,128-128,128C57.3,256,0,198.7,0,128C0,57.3,57.3,0,128,0C198.7,0,256,57.3,256,128z};
    \fill[white] svg{M86.3,186.2H70.9V79.1h15.4v48.4V186.2z}
                 svg{M108.9,79.1h41.6c39.6,0,57,28.3,57,53.6c0,27.5-21.5,53.6-56.8,53.6h-41.8V79.1z M124.3,172.4h24.5c34.9,0,42.9-26.5,42.9-39.7c0-21.5-13.7-39.7-43.7-39.7h-23.7V172.4z}
                 svg{M88.7,56.8c0,5.5-4.5,10.1-10.1,10.1c-5.6,0-10.1-4.6-10.1-10.1c0-5.6,4.5-10.1,10.1-10.1C84.2,46.7,88.7,51.3,88.7,56.8z};
  }
}

\newcommand{\SCORCIDScale}{0.00390625}

\newlength{\SCORCIDHeight}

\DeclareRobustCommand\ORCID[1]{%
\texorpdfstring{%
\setlength{\SCORCIDHeight}{\fontcharht\font`X}%
\href{https://orcid.org/#1}{\XeTeXLinkBox{\mbox{%
\begin{tikzpicture}[yscale=-\SCORCIDScale*\SCORCIDHeight,
xscale=\SCORCIDScale*\SCORCIDHeight,transform shape]
\pic{orcidlogo};
\end{tikzpicture}%
}}}}{}}

\graphicspath{{figure/}}
\setcitestyle{authoryear,round,aysep={},yysep={;}}
\theoremstyle{thmstyleone}
\newtheorem{theorem}{Theorem}

\newtheorem{lemma}{Lemma}
\theoremstyle{thmstyletwo}
\newtheorem{example}{Example}
\newtheorem{remark}{Remark}
\theoremstyle{thmstylethree}
\newtheorem{definition}{Definition}
\numberwithin{equation}{section}
\hypersetup{
  pdftitle={A unified spectral analysis of high-dimensional linear regression: debiasing and thresholding},
  pdfauthor={Ying-Ao Wang, Yunyi Zhang, Ye Zhang},
  pdfsubject={Statistics and Computing},
  hidelinks
}

\begin{document}
\title[Spectral analysis of high-dimensional regression]{A unified spectral analysis of high-dimensional linear regression: debiasing and thresholding}
\author[1]{\fnm{Ying-Ao} \sur{Wang}\,\ORCID{0000-0001-8073-380X}}\email{wya@bit.edu.cn}
\author[2]{\fnm{Yunyi} \sur{Zhang}\,\ORCID{0000-0001-7513-0908}}\email{zhangyunyi@cuhk.edu.cn}
\author*[1,3]{\fnm{Ye} \sur{Zhang}\,\ORCID{0000-0003-4023-6352}}\email{ye.zhang@smbu.edu.cn}
\affil[1]{\orgdiv{School of Mathematics and Statistics},
  \orgname{Beijing Institute of Technology},
  \orgaddress{\street{5 Zhongguancun South Street}, \city{Beijing}, \postcode{100081}, \country{China}}}
\affil[2]{\orgdiv{School of Data Science},
  \orgname{The Chinese University of Hong Kong, Shenzhen},
  \orgaddress{\street{2001 Longxiang Boulevard}, \city{Shenzhen}, \postcode{518172}, \country{China}}}
\affil[3]{\orgdiv{MSU-BIT-SMBU Joint Research Center of Applied Mathematics},
  \orgname{Shenzhen MSU-BIT University},
  \orgaddress{\street{1 International University Avenue}, \city{Shenzhen}, \postcode{518172}, \country{China}}}
\abstract{We use spectral analysis to establish a unified theory of debiasing, thresholding, and simultaneous inference for high-dimensional linear regression with ill-conditioned designs. The methods include Ridge regression, spectral cut-off, and iterative and dynamical regression procedures. A common debiasing procedure uses two applications of the original linear algorithm with fixed parameters. Spectral analysis shows that debiasing squares the residual, while implementations based on matrix-vector products require no explicit spectral decomposition. Under explicit design, signal, and moment conditions, we establish maximum-norm consistency, recovery of the population threshold set, and error bounds for coefficient and noise-variance estimation. We derive Gaussian approximations for maximum coordinate errors and, under full column rank and deterministic sequences of regression parameters and thresholds, prove bootstrap validity for simultaneous confidence regions. Under additional regularity and gain conditions, we relate noise-free bias rates to best-worst-case error rates. Sparse and non-sparse regression experiments examine the effects of debiasing and thresholding under specified tuning rules. An elliptic inverse source problem extends the numerical assessment to reconstruction and empirical simultaneous coverage.}
\maketitle
\section{Introduction}
\label{sec:Introduction}
We consider observations $(X_{i1},\cdots,X_{ip},Y_i)\in\mathbf{R}^p\times\mathbf{R}$, $i=1,\ldots,n$, satisfying the linear model
\begin{equation}
\label{LinearModel}
Y_i = \sum_{j = 1}^p X_{ij}\beta_j + e_i,
\end{equation}
where the errors $e_i$ are i.i.d. with $\mathbb{E}e_i=0$. When the dimension is comparable to the sample size and the design is ill-conditioned, small singular values can amplify stochastic error. Regularization controls this amplification at the cost of bias, while sparsity and subsequent inference introduce further requirements on the estimator and its tuning \citep{MR2807761,Statistics_data_science}. The statistical analysis must therefore account for regularization bias, random fluctuations, and the accuracy of any selected coordinate set.

This setting is closely connected to discretized inverse problems. Examples arise in computed tomography \citep{Hansen2021}, image deblurring \citep{8782833}, and inverse source reconstruction \citep{Zhang2020JCAM}. Their discretizations can produce large, ill-conditioned linear systems for which iterative methods offer alternatives to explicit matrix factorizations. A useful regression theory should accommodate these computational implementations while identifying conditions for consistent estimation and simultaneous inference.

We analyze a class of linear regression estimators using the generator-function representation from regularization theory \citep{engl1996regularization,Pereverzev2003,hofmann2007analysis,2019Onfractional}. Following its regression formulation in \citet{WANG2024101826}, we use the representation
\begin{equation}
\label{generalLRIntro}
  \frac{1}{n} g_{\alpha}(\frac{1}{n}\mathbf{X}_n^\mathrm{T}\mathbf{X}_n) \mathbf{X}_n^\mathrm{T}\mathbf{Y}_n,
\end{equation}
where $g_\alpha$ satisfies Definition~\ref{new} and $\alpha$ controls regularization. The associated residual $r_\alpha(\lambda)=1-\lambda g_\alpha(\lambda)$ describes the regularization bias. Bounds on these two functions provide a common basis for analyzing Ridge, spectral cut-off, and several iterative and dynamical methods. We develop the statistical analysis of their ordinary, debiased, and thresholded estimators as the dimension and conditioning of the design vary with sample size.

Spectral representation is a tool for analyzing the estimators; their computation follows the respective regression algorithms. The representation does not prescribe evaluation by spectral decomposition. For example, Landweber and Nesterov updates, as well as the explicit discretizations used for Showalter and HBF, can be evaluated through products with $\mathbf{X}_n$ and $\mathbf{X}_n^\mathrm{T}$. Stable step sizes may require an upper bound on the operator norm, but these updates do not require singular vectors or the full set of singular values. With the other numerical parameters fixed, the iteration count controls regularization. The stopping rule and numerical parameters remain part of the method specification.

The main contributions are as follows.
\begin{itemize}
    \item \textbf{A common statistical theory for linear regression estimators.} Under explicit conditions on the design, signal, noise, and regression parameter, we establish maximum-norm consistency, recovery of the population threshold set, Euclidean error bounds for thresholded coefficients, and consistency of noise-variance estimation (Theorems~\ref{thmP01}--\ref{them27}). The bounds distinguish regularization bias, stochastic error, discarded coefficients, and any unidentifiable component of the signal. This provides conditions that can be checked for different members of the class.
    \item \textbf{Unified debiasing and thresholding with a common implementation.} We develop a common debiasing procedure for the linear regression algorithms and combine it with coordinate thresholding. For a fixed linear algorithm, the debiased estimator can be computed by two applications of that algorithm with the same regression and numerical parameters; see \eqref{numachieve}. Thus implementations that avoid explicit spectral decomposition retain this property. Spectral analysis shows that debiasing replaces the residual $r_\alpha$ by $r_\alpha^2$ and separates the reduction in regularization bias from the change in the noise filter.
    \item \textbf{Gaussian approximation and simultaneous inference.} We establish Gaussian approximations for maximum coordinate errors of ordinary, debiased, and thresholded estimators (Theorems~\ref{them29}--\ref{them215}). Under full column rank and deterministic sequences of regression parameters and thresholds, Theorem~\ref{thm:bootstrap-validity} proves bootstrap validity for simultaneous confidence regions. The argument extends the Ridge bootstrap analysis of \citet{Zhang2020RidgeRR} to the class considered here.
\end{itemize}

Two further results support this analysis. First, under additional regularity and spectral gain conditions, we connect noise-free bias rates to best-worst-case error rates, using the rate transformation studied by \citet{boct2021convergence,WANG2024101826}. Second, controlled sparse and non-sparse regression experiments examine the effects of debiasing and thresholding under specified tuning rules, and an elliptic inverse source problem supplies an additional computational application.

Section~\ref{literature} relates the analysis to regularization theory and high-dimensional inference. Section~\ref{theory} defines the estimator class and its implementations, and Section~\ref{Reduction} develops the debiasing procedure. Section~\ref{Theory} gives the statistical and best-worst-case results. Section~\ref{simulation} presents the regression experiments and inverse source application, followed by conclusions in Section~\ref{sec:Con}. Full proofs are provided in Appendix~\ref{app:core-proofs} and the Supplementary Material (Online Resource~1), together with method-specific spectral derivations and numerical details.

\section{Related literature}
\label{literature}

We review two strands of literature that inform the methods and their statistical analysis: inverse problems and regularization theory, and high-dimensional estimation and statistical inference.

\paragraph*{Inverse problems and regularization theory.}
Regularization theory provides tools for stable approximation and convergence analysis in ill-posed problems \citep{tikhonov1963regularization,tikhonov1963solution,phillips1962technique,engl1996regularization,JinIto2014}. Generator functions, spectral residuals, and qualification bounds make it possible to compare regularization methods through common approximation properties \citep{Vainikko1986,Pereverzev2003,S0036142903420947,hofmann2007analysis}. Iterative regularization and discrepancy-based parameter choice are developed in \citet{morozov1966solution,Bakushinsky2004,kaltenbacher2008iterative,Hansen2010}; analyses of dynamical methods include \citet{2019Onfractional,boct2021convergence,YZhang2020On}. The regression framework of \citet{WANG2024101826} uses this generator formulation for several linear regression methods. These works supply the spectral tools and underlying algorithms used here. Our analysis connects their generator and residual bounds to estimation, thresholding, and simultaneous inference when the dimension grows with sample size.

\paragraph*{High-dimensional estimation and statistical inference.}
Regularization and variable selection have extensive statistical theory; see \citet{MR2807761,Statistics_data_science}. Lasso promotes sparse coefficient estimates through an $\ell_1$ penalty \citep{10.1111x,Knight2000AsymptoticsFL}. Model-selection theory includes nonconcave penalized likelihood \citep{MR1946581} and selection consistency for Lasso \citep{MR2274449}. Debiasing and inference for high-dimensional coefficients are developed in \citet{MR3102549,MR3153940}, while orthogonal-score constructions provide another approach to controlling the effect of nuisance estimation \citep{10.1111/ectj.12097}. Coordinate thresholding can be applied to fitted estimators, as in thresholded Lasso \citep{EJS624} and thresholded Ridge regression \citep{Zhang2020RidgeRR,MR4726969}. In our setting, linearity gives a common debiasing implementation using two applications of the original algorithm, and spectral analysis identifies the resulting bias and noise terms. The resulting analysis tracks both the regularization bias and the error contributed by discarded coefficients; recovery concerns the population threshold set under a coefficient-separation condition.

Gaussian approximation and bootstrap methods for maxima of high-dimensional sums provide tools for simultaneous inference \citep{13AOS1161,zhang2017simultaneous}. A particularly close antecedent is \citet{Zhang2020RidgeRR}, who study debiasing, thresholding, and bootstrap inference for Ridge regression. We extend this line of analysis to ordinary and debiased estimators described by common generator bounds. The Gaussian approximation results allow rank-deficient designs subject to the stated identifiability and variance conditions. The bootstrap validity theorem additionally requires full column rank and deterministic sequences of regression parameters and thresholds.

The contribution relative to these strands is a common construction and statistical analysis across a class of regularization methods, together with an implementation that uses the original linear algorithm. The method definitions and their cited qualification results provide the inputs to this analysis. The numerical comparisons then examine how debiasing and thresholding act under different signal structures and specified tuning rules, including an inverse source application.

\section{Debiasing and thresholding in linear regression}
\label{theory}

This manuscript focuses on the high-dimensional linear regression model
\begin{equation}
\label{hgfc}
\mathbf{Y}_n=\mathbf{X}_n\pmb{\beta}+\pmb{e}_n,
\end{equation}
where $\mathbf{Y}_n = [Y_1, \cdots, Y_n]^\mathrm{T} \in \mathbf{R}^n$ represents the response vector\footnote{The superscript ``$\mathrm{T}$'' means the transpose of a vector or a matrix.}, and $\mathbf{X}_n = [x_{ij}]_{1 \leq i \leq n, 1 \leq j \leq p} \in \mathbf{R}^{n \times p}$ denotes the fixed (non-random) design matrix. Its rank $s = s(n, p)$ may increase without bound as $n, p \to \infty$. $\pmb{\beta} = [\beta_1, \cdots,\beta_p]^\mathrm{T} \in \mathbf{R}^p$ denotes the coefficient vector to be estimated. $\pmb{e}_n = [e_1, \cdots, e_n]^\mathrm{T} \in \mathbf{R}^n$ denotes the error vector.

\subsection{A class of linear regression methods}

We use the generator--residual formulation of general linear regularization \citep{Pereverzev2003,hofmann2007analysis,2019Onfractional}, following its regression formulation in \citet{WANG2024101826}. The pair $(g_\alpha,r_\alpha)$ describes the estimator and its regularization bias, while qualification bounds control the decay of the residual as the regression parameter decreases. These standard tools provide the method-specific inputs to the common statistical analysis in Section~\ref{Theory}.

\begin{definition}
\label{new}
  A family of functions $g_{\alpha}(\lambda)$ $(\lambda>0)$, defined for regression parameters $0\le \alpha \le \overline{\alpha}$,
  constitutes a generator of linear regression methods for problem \eqref{hgfc} if the following three conditions are satisfied:
  \begin{enumerate}[(D1-3)]
    \item[(D1-1)] For the bias function $r_{\alpha}(\lambda):=1- \lambda g_{\alpha}(\lambda)$, it holds that $\lim\limits_{\alpha\to 0}|r_{\alpha}(\lambda)|=0$ for any fixed $\lambda\in (0,+\infty]$.
    \item[(D1-2)] There exists a constant $c_r> 0$ such that $|r_{\alpha}(\lambda)|\le c_r$ for all $\lambda\in (0,+\infty]$ and $\alpha\in (0,\overline{\alpha}]$.
    \item[(D1-3)] There exists a constant $c_0> 0$ such that $ g_{\alpha}(\lambda)\leq \min (2/\lambda, c_0/\sqrt{\lambda \alpha})$ for all $\lambda\in (0,+\infty]$ and $\alpha\in (0,\overline{\alpha}]$.
  \end{enumerate}
\end{definition} 
The linear regression estimators considered here admit the common representation
\begin{equation}
\label{generalLR}
  \hat{\pmb{\beta}}_{\alpha}=\frac{1}{n}g_{\alpha}(\frac{1}{n}\mathbf{X}_n^\mathrm{T}\mathbf{X}_n)\mathbf{X}_n^\mathrm{T}\mathbf{Y}_n.
\end{equation}

Regularization bias can be substantial relative to stochastic error in high-dimensional regression and can affect the coverage of confidence regions \citep{MR3153940,10.1111/ectj.12097,Zhang2020RidgeRR}. Its magnitude depends on the design, the signal, and the regression parameter. Section~\ref{Reduction} derives a common debiasing procedure for the present class of estimators. The resulting debiased estimator, given below and in \eqref{debiasLR}, is evaluated together with its stochastic error rather than through bias reduction alone.
\begin{equation*}
  \tilde{\pmb{\beta}}_{\alpha}
= \frac{1}{n}\left[ \mathbf{I}+ r_{\alpha}(\frac{1}{n}\mathbf{X}_n^\mathrm{T}\mathbf{X}_n) \right]g_{\alpha}(\frac{1}{n}\mathbf{X}_n^\mathrm{T}\mathbf{X}_n)\mathbf{X}_n^\mathrm{T}\mathbf{Y}_n.
\end{equation*}

\begin{definition}[Slightly modified from Definition 2.3 in \citet{YZhang2020On}]
\label{indexd}
  A linear regression method \eqref{generalLR} for model \eqref{hgfc}, with generator $g_{\alpha}(\lambda)$ for $0<\lambda\leq C_{\lambda}$,
  is said to have monomial qualification of order $d$ if the following inequality holds
\begin{equation}
\label{indd}
    \sup _{\lambda \in\left(0, C_{\lambda}\right]}\left|r_\alpha(\lambda)\right| \lambda^d \leq C_* \alpha^d,
  \end{equation}
where $C_{\lambda}$ and $C_*$ are constants independent of the value of $\alpha$.
\end{definition}

\begin{remark}
Definitions~\ref{new} and~\ref{indexd} are standard tools in regularization theory \citep{engl1996regularization}. Here they provide common bounds for separating regularization bias from stochastic error. Combined with conditions on the design, signal, noise, and regression parameter, these bounds yield the estimation and inference results in Section~\ref{Theory}.
\end{remark}

For a positive threshold $b_n$, which may depend on the sample size, let $\mathcal{N}_{b_n}$ denote the set of indices $i$ for which $|\beta_i|>b_n$, i.e., $\mathcal{N}_{b_n}=\left\{i || \beta_i |>b_n\right\}$. Similarly, we define the index sets $\widehat{\mathcal{N}}_{b_n}$ and $\widetilde{\mathcal{N}}_{b_n}$, the thresholded estimator \mbox{$\hat{\pmb{\theta}}=\left(\hat{\theta}_1, \cdots, \hat{\theta}_p\right)^\mathrm{T}$} and the thresholded debiased estimator \mbox{$\tilde{\pmb{\theta}}=\left(\tilde{\theta}_1, \cdots, \tilde{\theta}_p\right)^\mathrm{T}$} as follows:
\begin{equation}
\label{threshold-beta}
\begin{aligned}
\widehat{\mathcal{N}}_{b_n}&=\left\{i|| (\hat{\pmb{\beta}}_\alpha)_i |>b_n\right\},\\
\hat{\theta}_i&=(\hat{\pmb{\beta}}_\alpha)_i\times\pmb{\mathbf{1}}_{i\in\widehat{\mathcal{N}}_{b_n}},\\
\widetilde{\mathcal{N}}_{b_n}&=\left\{i|| (\tilde{\pmb{\beta}}_\alpha)_i |>b_n\right\},\\
\tilde{\theta}_i&=(\tilde{\pmb{\beta}}_\alpha)_i\times\pmb{\mathbf{1}}_{i\in\widetilde{\mathcal{N}}_{b_n}}.
\end{aligned}
\end{equation}
Coordinate thresholding sets small fitted coefficients to zero. For sparse signals, this can reduce the accumulation of estimation error over inactive coordinates. The resulting Euclidean error also depends on the true coefficients discarded by thresholding. Theorem~\ref{Thm1} bounds both contributions, while Theorem~\ref{thmP1} gives recovery of the population threshold set under coefficient separation.

Simultaneous confidence regions require control of the maximum estimation error across coordinates. We use rectangular regions based on $\|\hat{\pmb\theta}-\pmb\beta\|_\infty$ and $\|\tilde{\pmb\theta}-\pmb\beta\|_\infty$, whose quantiles are approximated by bootstrap simulation. Algorithm~\ref{algorithm.bootstrap} follows the resampling approach of \citet{Zhang2020RidgeRR}, with the ordinary and debiased generators of the present class. Gaussian approximation provides the theoretical link between the original and bootstrap statistics \citep{13AOS1161,zhang2017simultaneous}.

\begin{algorithm*}[t]
\caption{Wild bootstrap}\label{algorithm.bootstrap}
\begin{minipage}{\linewidth}
\textbf{Input}: Design matrix $\mathbf{X}_n$, response vector $\mathbf{Y}_n=\mathbf{X}_n\pmb{\beta}+\pmb{e}_n$, threshold $b_n$\footnote{The numerical section specifies the original-sample threshold selection for each method and estimator branch.}, nominal coverage probability $1-\alpha^*$, number of bootstrap replicates $B$
    \begin{enumerate}
        \item Calculate $\hat{\pmb{\theta}}$ and $\tilde{\pmb{\theta}}$\footnote{Original-sample tuning follows the experiment-specific rules stated in the numerical section.} defined in \eqref{threshold-beta}, along with 
        \begin{equation}
\label{threshold-sigma}
\widehat{\sigma}^2=\frac{1}{n} \sum_{i=1}^n\left(Y_i-\sum_{j=1}^p x_{i j} \hat{\theta}_j\right)^2 \quad \text{and} \quad  \widetilde{\sigma}^2=\frac{1}{n} \sum_{i=1}^n\left(Y_i-\sum_{j=1}^p x_{i j} \tilde{\theta}_j\right)^2.
        \end{equation}
        \item Generate the error vector $\hat{\pmb{e}}_n=[\hat{e}_1, \cdots, \hat{e}_n]^\mathrm{T}$ whose components are i.i.d. $N(0,\widehat{\sigma}^2)$ variables. Then calculate $\hat{\mathbf{Y}}^*_n=\mathbf{X}_n\hat{\pmb{\theta}}+\hat{\pmb{e}}_n$. Similarly, $\tilde{\pmb{e}}_n$ and $\tilde{\mathbf{Y}}^*_n$ are generated following the same process.
        \item Hold the regression parameter $\alpha$ selected on the original sample fixed throughout the bootstrap replicates; for a discrete iterative method, hold its selected iteration count and numerical parameters fixed. Calculate $\hat{\pmb{\beta}}^*_{\alpha}=\frac{1}{n}g_{\alpha}(\frac{1}{n}\mathbf{X}_n^\mathrm{T}\mathbf{X}_n)\mathbf{X}_n^\mathrm{T}\hat{\mathbf{Y}}^*_n$ and $\tilde{\pmb{\beta}}^*_{\alpha}=\frac{1}{n}[ \mathbf{I}+  r_{\alpha}(\frac{1}{n}\mathbf{X}_n^\mathrm{T}\mathbf{X}_n)]g_{\alpha}(\frac{1}{n}\mathbf{X}_n^\mathrm{T}\mathbf{X}_n)\mathbf{X}_n^\mathrm{T}\tilde{\mathbf{Y}}^*_n$. The threshold $b_n$ is also held fixed, but the retained coordinate sets are recalculated in the next step.
        \item Calculate $\widehat{\mathcal{N}}^*_{b_n}=\left\{i|| (\hat{\pmb{\beta}}^*_\alpha)_i |>b_n\right\}$ and $\hat{\pmb{\theta}}^*=\left[\hat{\theta}^*_1, \cdots, \hat{\theta}^*_p\right]^\mathrm{T}$ with $\hat{\theta}^*_i=(\hat{\pmb{\beta}}^*_\alpha)_i \times \pmb{\mathbf{1}}_{i \in \widehat{\mathcal{N}}^*_{b_n}}$ for $i=1, \cdots, p$. Define $\widetilde{\mathcal{N}}^*_{b_n}$ and $\tilde{\pmb{\theta}}^*$ analogously.
        \item Calculate $\hat{E}_b^*= \max\limits_{i=1, \cdots, p}\left|\hat{\theta}^*_i-\hat{\theta}_i\right|$ and $\tilde{E}_b^*= \max\limits_{i=1, \cdots, p}\left|\tilde{\theta}^*_i-\tilde{\theta}_i\right|$.
        \item (For constructing a confidence region) Repeat steps 2 to 5 $B$ times to generate $\hat{E}_b^*, \tilde{E}_b^*, b=1,2, \ldots, B$; then calculate the $1-\alpha^*$ sample quantile $\hat{C}_{1-\alpha^*}^*$ of $\hat{E}_b^*$ and $\tilde{C}_{1-\alpha^*}^*$ of $\tilde{E}_b^*$. The $1-\alpha^*$ confidence regions for the parameter of interest $\pmb{\beta}$ are given by the sets
        \begin{equation}
\label{CR-1}
\begin{gathered}
\left\{\pmb{\beta}=(\beta_1,\cdots,\beta_p)^\mathrm{T}\vert
\max_{i=1,\cdots,p}|\hat{\theta}_i-\beta_i|\leq\hat{C}_{1-\alpha^*}^*\right\}
\text{ and }
\left\{\pmb{\beta}=(\beta_1,\cdots,\beta_p)^\mathrm{T}\vert
\max_{i=1,\cdots,p}|\tilde{\theta}_i-\beta_i|\leq\tilde{C}_{1-\alpha^*}^*\right\}.
\end{gathered}
\end{equation}
    \end{enumerate}

\end{minipage}
\end{algorithm*}
Under its stated assumptions, including full column rank and deterministic sequences of regression parameters and thresholds, Theorem~\ref{thm:bootstrap-validity} gives an asymptotic coverage guarantee for Algorithm~\ref{algorithm.bootstrap}.
\subsection{Examples of linear regression methods}
\label{example}
We now define ten regression methods and describe their debiased and thresholded counterparts within the preceding framework. The underlying methods have established origins: the spectral representations of least squares, Ridge, Landweber, Showalter, second-order asymptotic regression (SOAR), and heavy ball with friction (HBF) regression were presented in \citet[Examples 3.1--3.6]{WANG2024101826}. The sources for spectral cut-off (SC), fractional asymptotical regression (FAR), acceleration regression of order $\kappa$ (AR$^\kappa$), and discrete Nesterov acceleration are cited in the corresponding examples below. Our construction applies the common debiasing and thresholding operations to these methods; the discrete Nesterov iteration is distinguished from the continuous SOAR flow.

Write $\mathbf{A}_n=\mathbf{X}_n^\mathrm{T}\mathbf{X}_n/n$ and $\pmb{c}_n=\mathbf{X}_n^\mathrm{T}\mathbf{Y}_n/n$. For $\lambda>0$, each residual function determines both the ordinary and debiased generator functions:
\begin{equation}
\label{method-filters}
\begin{aligned}
g_\alpha(\lambda)&=\frac{1-r_\alpha(\lambda)}{\lambda},\\
G_\alpha(\lambda)&=[1+r_\alpha(\lambda)]g_\alpha(\lambda)=\frac{1-r_\alpha^2(\lambda)}{\lambda}.
\end{aligned}
\end{equation}
Thus $\hat{\pmb{\beta}}_\alpha=g_\alpha(\mathbf{A}_n)\pmb{c}_n$ and $\tilde{\pmb{\beta}}_\alpha=G_\alpha(\mathbf{A}_n)\pmb{c}_n$. For every example below, coordinate thresholding of these two estimators is defined by \eqref{threshold-beta}. Table~\ref{tab:method-summary} summarizes their residual functions and qualification ranges. Method-specific spectral derivations and qualification bounds are given in Section~\ref{supp-supp:sec:filters} of the Supplementary Material (Online Resource~1); detailed numerical discretizations are given in Section~\ref{supp-supp:sec:implementation}.

All method definitions below use the normalized quantities $\mathbf{A}_n$ and $\pmb{c}_n$. The numerical schemes in the Supplementary Material are written for a generic design--response pair: supplying $(\mathbf{X}_n/\sqrt n,\mathbf{Y}_n/\sqrt n)$ gives the normalized implementation. Section~\ref{simulation} describes the stopping and parameter-selection rules; the numerical settings for the unscaled experimental inputs are given in Supplementary Section~\ref{supp-supp:sec:regression-settings}.

\begin{table*}[t]
\centering
\caption{Residual functions and verified qualification orders for the methods in Section~\ref{example}.}
\label{tab:method-summary}

\begingroup
\small
\setlength{\tabcolsep}{4pt}
\renewcommand{\arraystretch}{1.25}
\begin{tabular}{p{0.17\linewidth}p{0.38\linewidth}p{0.22\linewidth}p{0.11\linewidth}}
\hline
Method & Residual $r_\alpha(\lambda)$ & Verified orders $d$ & Example\\
\hline
LS (baseline) & $0$ on the positive spectrum & --- & \ref{ex:LS}\\
SC & $\mathbf{1}_{\{\lambda<\alpha\}}$ & Any $d>0$ & \ref{ex:SC}\\
Ridge & $\alpha/(\lambda+\alpha)$ & $0<d\leq1$ & \ref{ex:ridge-summary}\\
Landweber & $(1-h\lambda)^k$, $k=\lfloor1/\alpha\rfloor$ & Any $d>0$ & \ref{ex3}\\
Showalter & $e^{-\lambda/\alpha}$ & Any $d>0$ & \ref{ex:Showalter}\\
SOAR & $2^{s^*}\Gamma(s^*+1)J_{s^*}(z)/z^{s^*}$ & $0<d\leq(1+2s^*)/4$ & \ref{ex:SOAR}\\
HBF & $R_\eta(1/\alpha,\lambda)$ & Any $d>0$ & \ref{ex:HBF}\\
FAR & $E_\vartheta(-\lambda/\alpha)$ & $0<d\leq1$ & \ref{ex:FAR}\\
AR$^\kappa$ & $e^{-\lambda/\alpha}$ & Any $d>0$ & \ref{ex:ARk}\\
Nesterov & $r_k(\lambda)$ in \eqref{Nesterov-residual-main}, $\alpha=k^{-2}$ & $0<d\leq(\omega+1)/4$ & \ref{ex:Nesterov}\\
\hline
\end{tabular}
\par\smallskip
\begin{minipage}{0.97\linewidth}
\footnotesize
\textit{Notes.} Here $\lambda\in(0,C_\lambda]$; notation and parameter choices are specified in the cited examples.
The listed residuals and qualification orders refer to the ordinary estimators; debiasing replaces $r_\alpha$ by $r_\alpha^2$.
The orders refer to the residual bound~\eqref{indd}, not the additional regularity conditions of Definition~\ref{DefRegular}. For FAR, the listed orders apply to each fixed $\vartheta\in(0,1)\cup(1,2)$; at $\vartheta=1$, it coincides with Showalter regression and has qualification of every order $d>0$.
\end{minipage}
\endgroup
\end{table*}

\begin{example}[Least squares]
\label{ex:LS}
The minimum-norm least-squares estimator is
\begin{equation}
\label{LS-main}
\hat{\pmb{\beta}}_{\mathrm{LS}}=\mathbf{A}_n^\dagger\pmb{c}_n,
\end{equation}
where $\dagger$ denotes the Moore--Penrose inverse. For full column rank this is the usual least-squares solution. On the positive spectrum, $g(\lambda)=1/\lambda$ and $r(\lambda)=0$, so debiasing leaves this estimator unchanged. Coordinate thresholding can still be applied. LS is included as a baseline: its generator does not satisfy the uniform bound (D1-3), and the zero residual alone does not establish the assumptions of the statistical theorems.
\end{example}

\begin{example}[Spectral cut-off]
\label{ex:SC}
For a cut-off $\alpha>0$, spectral cut-off regression \citep[Example 5]{BAUER200752} uses
\begin{equation}
\label{SC-main}
g_\alpha(\lambda)=\frac{\mathbf{1}_{\{\lambda\geq\alpha\}}}{\lambda},
\qquad r_\alpha(\lambda)=\mathbf{1}_{\{\lambda<\alpha\}}.
\end{equation}
The estimator retains only the eigendirections of $\mathbf{A}_n$ with eigenvalues at least $\alpha$. Since $r_\alpha^2=r_\alpha$, we have $G_\alpha=g_\alpha$: at a fixed cut-off, debiasing leaves the unthresholded estimator unchanged. Spectral truncation and coordinate thresholding are distinct operations. SC has qualification of every order $d>0$, but it does not satisfy the additional regularity conditions of Definition~\ref{DefRegular}; the verification is given in Supplementary Example~\ref{supp-method:SC}.
\end{example}

\begin{example}[Ridge regression]
\label{ex:ridge-summary}
For $g_\alpha(\lambda)=(\lambda+\alpha)^{-1}$, the normalized formulation gives
\begin{equation}
\label{ridge-summary}
\begin{aligned}
\hat{\pmb{\beta}}_\alpha&=(\mathbf{A}_n+\alpha\mathbf{I})^{-1}\pmb{c}_n,\\
\tilde{\pmb{\beta}}_\alpha&=(\mathbf{A}_n+2\alpha\mathbf{I})
(\mathbf{A}_n+\alpha\mathbf{I})^{-2}\pmb{c}_n.
\end{aligned}
\end{equation}
The first expression minimizes
$n^{-1}\|\mathbf{Y}_n-\mathbf{X}_n\pmb{\beta}\|^2+\alpha\|\pmb{\beta}\|^2$.
Consequently, $\alpha$ here corresponds to the penalty $n\alpha$ when the squared loss is not divided by $n$. Its qualification range is $0<d\leq1$; the scalar bounds are given in Example~\ref{supp-method:Ridge}. Thresholding either estimator uses \eqref{threshold-beta}.
\end{example}

\begin{example}[Landweber regression]
\label{ex3}
With normalized inputs, the Landweber iteration of \citet[Example 3.3]{WANG2024101826} reads
\begin{equation}
\label{LandweberEX}
\begin{cases}
\pmb{\beta}_{k+1}=\pmb{\beta}_k+h(\pmb{c}_n-\mathbf{A}_n\pmb{\beta}_k),\\
\pmb{\beta}_0=\pmb{0},\qquad 0<h\leq1/C_\lambda.
\end{cases}
\end{equation}
Its residual is $r_k(\lambda)=(1-h\lambda)^k$, with $g_k(\lambda)=[1-r_k(\lambda)]/\lambda$. Since
$[1+r_k(\lambda)]g_k(\lambda)=g_{2k}(\lambda)$, the debiased estimator is obtained by doubling the number of iterations:
\begin{equation}
\label{debiasedLandweber}
\begin{gathered}
\hat{\pmb{\beta}}_k=\pmb{\beta}_k,\qquad
\tilde{\pmb{\beta}}_k=\pmb{\beta}_{2k},\\
\tilde{\theta}_i=(\pmb{\beta}_{2k})_i
\mathbf{1}_{\{i\in\widetilde{\mathcal{N}}_{b_n}\}}.
\end{gathered}
\end{equation}
Both runs use the same step size. For an implementation written with $\mathbf{X}_n^\mathrm{T}(\mathbf{Y}_n-\mathbf{X}_n\pmb{\beta}_k)$, the corresponding step size is $\Delta t=h/n$.
\end{example}

\begin{example}[Showalter regression]
\label{ex:Showalter}
The continuous-time counterpart of Landweber regression is the Showalter flow \citep[Example 3.4]{WANG2024101826}:
\begin{equation}
\label{first}
\begin{cases}
\dot{\pmb{\beta}}(t)+\mathbf{A}_n\pmb{\beta}(t)=\pmb{c}_n,
\qquad t>0,\\
\pmb{\beta}(0)=\pmb{0}.
\end{cases}
\end{equation}
At time $t=1/\alpha$, its residual is $r_\alpha(\lambda)=e^{-\lambda/\alpha}$. Thus $G_\alpha=g_{\alpha/2}$ and the debiased estimator is the exact flow evaluated at time $2t$. This identity concerns the continuous flow; a discretized implementation uses fixed-parameter debiasing as in \eqref{numachieve}. We use a fourth-order Runge--Kutta scheme, detailed in Supplementary Section~\ref{supp-supp:sec:implementation}.
\end{example}

\begin{example}[Second-order asymptotic regression]
\label{ex:SOAR}
SOAR uses a second-order flow with vanishing damping \citep{gong2020new,WANG2024101826}:
\begin{equation}
\label{second}
\begin{cases}
\ddot{\pmb{\beta}}(t)+\frac{1+2s^*}{t}\dot{\pmb{\beta}}(t)
+\mathbf{A}_n\pmb{\beta}(t)=\pmb{c}_n,\\
\pmb{\beta}(0)=\pmb{0},\qquad
\dot{\pmb{\beta}}(0)=\pmb{0},
\end{cases}
\end{equation}
where $t>0$ and $s^*>-1/2$ is fixed. Set $\alpha=\rho/t^2$ for a fixed $\rho>0$ and write $z=\sqrt{\rho\lambda/\alpha}$. The residual in Table~\ref{tab:method-summary} involves the Bessel function $J_{s^*}$ of the first kind; its value at zero is understood by continuity. The debiased estimator uses $G_\alpha=(1-r_\alpha^2)/\lambda$, followed by \eqref{threshold-beta}. The qualification range is $0<d\leq(1+2s^*)/4$. The discrete update used in the regression experiments is recorded in Supplementary Section~\ref{supp-supp:sec:implementation}.
\end{example}

\begin{example}[Heavy ball with friction regression]
\label{ex:HBF}
HBF replaces the vanishing damping in SOAR by a fixed coefficient $\eta>0$ \citep{WANG2024101826,YZhang2020On}:
\begin{equation}
\label{ball}
\begin{cases}
\ddot{\pmb{\beta}}(t)+\eta\dot{\pmb{\beta}}(t)
+\mathbf{A}_n\pmb{\beta}(t)=\pmb{c}_n,\\
\pmb{\beta}(0)=\pmb{0},\qquad
\dot{\pmb{\beta}}(0)=\pmb{0}.
\end{cases}
\end{equation}
With $t=1/\alpha$, the residual is $R_\eta(t,\lambda)$, which solves $\partial_t^2R_\eta+\eta\partial_tR_\eta+\lambda R_\eta=0$ with $R_\eta(0,\lambda)=1$ and $\partial_tR_\eta(0,\lambda)=0$. Its overdamped, critically damped, and underdamped expressions are given in Supplementary Example~\ref{supp-method:HBF}. The debiased generator is $G_\alpha=(1-R_\eta^2(1/\alpha,\lambda))/\lambda$. For fixed $\eta$ and bounded spectrum, HBF has qualification of every order $d>0$. Its numerical implementation applies a fourth-order Runge--Kutta scheme to the equivalent first-order system.
\end{example}

\begin{example}[Fractional asymptotical regression]
\label{ex:FAR}
FAR replaces the first derivative in the Showalter flow by a Caputo derivative of order $\vartheta$ \citep{2019Onfractional}:
\begin{equation}
\label{frac}
\begin{cases}
({}^C D_{0+}^{\vartheta}\pmb{\beta})(t)
+\mathbf{A}_n\pmb{\beta}(t)=\pmb{c}_n,\\
\pmb{\beta}^{(j)}(0)=\pmb{0},\qquad j=0,\ldots,m-1,
\end{cases}
\end{equation}
where $\vartheta\in(0,1)\cup(1,2)$ and $m=\lceil\vartheta\rceil$. Here
\begin{equation*}
({}^C D_{0+}^{\vartheta}\pmb{\beta})(t)
=\frac{\displaystyle\int_0^t(t-u)^{m-\vartheta-1}
\pmb{\beta}^{(m)}(u)\,du}{\Gamma(m-\vartheta)}.
\end{equation*}
For $\vartheta=1$, this is the ordinary first derivative with the sole initial condition $\pmb{\beta}(0)=\pmb{0}$, recovering Showalter regression. Setting $\alpha=t^{-\vartheta}$ gives
\begin{equation}
\label{FAR-residual-main}
\begin{aligned}
r_\alpha(\lambda)&=E_\vartheta(-\lambda/\alpha),\\
G_\alpha(\lambda)&=\frac{1-E_\vartheta^2(-\lambda/\alpha)}{\lambda},
\end{aligned}
\end{equation}
where $E_\vartheta(z)=\sum_{j=0}^{\infty}z^j/\Gamma(\vartheta j+1)$ is the Mittag--Leffler function. For each fixed $\vartheta\in(0,1)\cup(1,2)$, FAR has monomial qualification $0<d\leq1$, with maximal order $1$; for $\vartheta=1$, every order $d>0$ is admissible. Supplementary Example~\ref{supp-method:FAR} gives the generator and qualification bounds and separately verifies the additional conditions of Definition~\ref{DefRegular} for $0<\vartheta<1$. The fractional predictor--corrector scheme and its quadrature weights are given in Supplementary Section~\ref{supp-supp:sec:implementation}.
\end{example}

\begin{example}[Acceleration regression of order $\kappa$]
\label{ex:ARk}
Writing $\pmb{\beta}=\pmb{\beta}(t)$, the second-order dynamics associated with AR$^\kappa$ \citep{zhang2023acceleration} read
\begin{equation}
\label{orderk}
t\ddot{\pmb{\beta}}+(t^{-\kappa}-\kappa)\dot{\pmb{\beta}}
+\mathbf{A}_n\bigl[\pmb{\beta}+t^{\kappa+1}\dot{\pmb{\beta}}\bigr]=\pmb{c}_n,
\end{equation}
for $t>0$ and $\kappa>-1$. We use the solution branch characterized by
\begin{equation*}
\begin{cases}
\dot{\pmb{\beta}}(t)=t^\kappa(\pmb{c}_n-\mathbf{A}_n\pmb{\beta}(t)),\\
\pmb{\beta}(0)=\pmb{0}.
\end{cases}
\end{equation*}
This branch also satisfies $\dot{\pmb{\beta}}(0)=\pmb{0}$ when $\kappa>0$. Its residual is $\exp\{-\lambda t^{\kappa+1}/(\kappa+1)\}$. Consequently, with $\alpha=(\kappa+1)/t^{\kappa+1}$, its generator is the Showalter generator under a change of time. In the exact flow, debiasing therefore replaces $t$ by $2^{1/(\kappa+1)}t$. The numerical implementation of AR$^\kappa$ used in the regression experiments is given in Supplementary Section~\ref{supp-supp:sec:implementation}.
\end{example}

\begin{example}[Nesterov acceleration regression]
\label{ex:Nesterov}
The discrete Nesterov method combines an extrapolation step with a gradient step \citep{neubauer2017nesterov,Kindermann2021}. A related discretization was used in \citet[Section 4.1, equation (55)]{WANG2024101826}. For $\omega>-1$ and fixed $0<h\leq1/C_\lambda$, define
\begin{equation}
\label{Nesterov}
\begin{cases}
\pmb{z}_k=\pmb{\beta}_k+\frac{k-1}{k+\omega}
(\pmb{\beta}_k-\pmb{\beta}_{k-1}),\\
\pmb{\beta}_{k+1}=\pmb{z}_k+h(\pmb{c}_n-\mathbf{A}_n\pmb{z}_k),
\qquad k\geq1,\\
\pmb{\beta}_0=\pmb{0},\qquad \pmb{\beta}_1=h\pmb{c}_n.
\end{cases}
\end{equation}
Its residual polynomials are determined by
\begin{equation}
\label{Nesterov-residual-main}
\begin{cases}
r_0(\lambda)=1,\qquad r_1(\lambda)=1-h\lambda,\\[2pt]
\begin{aligned}[t]
r_{k+1}(\lambda)={}&
\frac{(2k+\omega-1)(1-h\lambda)}{k+\omega}\,r_k(\lambda)\\
&-\frac{(k-1)(1-h\lambda)}{k+\omega}\,r_{k-1}(\lambda).
\end{aligned}
\end{cases}
\end{equation}
On the parameter grid $\alpha=k^{-2}$, its qualification range is $0<d\leq(\omega+1)/4$. The corresponding Gegenbauer-polynomial representation and qualification analysis are given in Supplementary Example~\ref{supp-method:Nesterov}. The ordinary estimator is $\pmb{\beta}_k$; its debiased version uses \eqref{numachieve} with the same $k$, $h$, and $\omega$, followed by coordinate thresholding. For unscaled inputs, the step size in the gradient update is $\Delta t=h/n$.
\end{example}

\noindent\textbf{Common implementation of debiasing.}
All the linear iterations considered here admit a common implementation of debiasing as specified in \eqref{method-filters}. Let $\pmb{\beta}_k(\mathbf{X}_n,\mathbf{y})$ denote the output of a fixed linear algorithm after $k$ steps, with a zero initial state and fixed numerical parameters. Then
\begin{equation}
\label{numachieve}
\begin{aligned}
\hat{\pmb{\beta}}_k&=\pmb{\beta}_k(\mathbf{X}_n,\mathbf{Y}_n),\\
\tilde{\pmb{\beta}}_k&=2\hat{\pmb{\beta}}_k
-\pmb{\beta}_k(\mathbf{X}_n,\mathbf{X}_n\hat{\pmb{\beta}}_k).
\end{aligned}
\end{equation}
The second call uses the same iteration count, step size, and other algorithmic parameters as the first; its stopping rule is not reselected using the fitted response. This computes the debiased estimator without forming a spectral matrix function. If the underlying algorithm is implemented without spectral decomposition, the second call has the same property. The stopping rules used to choose the first run are described in Section~\ref{simulation}, while method-specific discretizations are supplied in Section~\ref{supp-supp:sec:implementation} of the Supplementary Material.

\section{Debiasing}
\label{Reduction}
In this section, we present the construction of the debiased regression estimator. To achieve this, we first define the noise-free intermediate quantity $\pmb{\beta}_{\alpha}:=\frac{1}{n}g_{\alpha}(\mathbf{A}_n)\mathbf{X}_n^\mathrm{T}\mathbf{X}_n\pmb{\beta}$. 
Then, the total error of linear regression \eqref{generalLR} can be decomposed as:
\begin{equation*}
    \begin{aligned}
      \hat{\pmb{\beta}}_{\alpha}-\pmb{\beta}  &=\hat{\pmb{\beta}}_{\alpha}-\pmb{\beta}_{\alpha}+\pmb{\beta}_{\alpha}-\pmb{\beta}\\
      &=\frac{1}{n}g_{\alpha}(\mathbf{A}_n)\mathbf{X}_n^\mathrm{T}\pmb{e}_n -  r_{\alpha}(\mathbf{A}_n)\pmb{\beta}. 
    \end{aligned}
\end{equation*}
The bias of $\hat{\pmb{\beta}}_{\alpha}$ is $-r_{\alpha}(\mathbf{A}_n)\pmb{\beta}$. We estimate this bias by $-r_{\alpha}(\mathbf{A}_n)\hat{\pmb{\beta}}_{\alpha}$.
Subtracting this estimate from $\hat{\pmb{\beta}}_{\alpha}$ gives the debiased estimator:
\begin{equation}
\label{debiasLR}
\begin{aligned}
\tilde{\pmb{\beta}}_{\alpha}
&=\hat{\pmb{\beta}}_{\alpha}+r_{\alpha}(\mathbf{A}_n)\hat{\pmb{\beta}}_{\alpha}\\
&=\frac{1}{n}[\mathbf{I}+r_{\alpha}(\mathbf{A}_n)]
 g_{\alpha}(\mathbf{A}_n)\mathbf{X}_n^\mathrm{T}\mathbf{Y}_n.
\end{aligned}
\end{equation}

Using \eqref{debiasLR}, the error decomposition above, and the generator function
$G_\alpha$ in \eqref{method-filters}, we obtain

\begin{equation}
\label{debiasest}
\begin{aligned}
\tilde{\pmb{\beta}}_{\alpha}-\pmb{\beta}
&=(\hat{\pmb{\beta}}_{\alpha}-\pmb{\beta})
+r_{\alpha}(\mathbf{A}_n)\hat{\pmb{\beta}}_{\alpha}\\
&=\frac{1}{n}G_{\alpha}(\mathbf{A}_n)\mathbf{X}_n^\mathrm{T}\pmb{e}_n
-r^2_{\alpha}(\mathbf{A}_n)\pmb{\beta}.
\end{aligned}
\end{equation}

For a deterministic choice of $\alpha$, if $|r_{\alpha}(\cdot)|\leq 1$, then
\begin{equation*}
       \begin{aligned}
\left\lVert \mathbb{E}  \tilde{\pmb{\beta}}_{\alpha}-\pmb{\beta}\right\rVert_2 &= \left\lVert r^2_{\alpha}(\mathbf{A}_n)\pmb{\beta}\right\rVert_2\\
&\leq \left\lVert r_{\alpha}(\mathbf{A}_n)\pmb{\beta}\right\rVert_2\\
&=\left\lVert \mathbb{E}  \hat{\pmb{\beta}}_{\alpha}-\pmb{\beta}\right\rVert_2.
\end{aligned}
\end{equation*}

On the identifiable spectral subspace, debiasing replaces $r_\alpha$ by $r_\alpha^2$ in the bias term and multiplies the noise filter by $1+r_\alpha$. The qualification and parameter conditions in Section~\ref{Theory} control these two effects, allowing a smaller order of regularization bias while retaining the same order of stochastic error. Finite-sample variance is governed by the transformed noise filter, and mean squared error combines this variance with squared bias. The identifiability conditions in that section control the component of $\pmb{\beta}$ in the null space of the design, which remains unchanged by debiasing. The numerical experiments assess how these effects interact with the signal, regularization choice, and coordinate thresholding.

\begin{remark}
   Debiasing has been developed through several constructions, including node-wise regression \citep{MR3153940} and orthogonal scores \citep{10.1111/ectj.12097}. Debiasing in \eqref{debiasLR} uses the linearity of the regression algorithm. Its spectral representation is used to analyze bias and stochastic error, while \eqref{numachieve} implements it through two applications of the original algorithm with the same regression and numerical parameters. Section~\ref{Theory} establishes statistical guarantees for this construction under the stated assumptions.
\end{remark}

\section{Theoretical results}
\label{Theory}

We first establish estimation and threshold-set consistency, then derive Gaussian approximations and bootstrap coverage guarantees. Proof sketches accompany Theorems~\ref{thmP01}, \ref{them29}, \ref{them215}, and~\ref{thm:bootstrap-validity}. Appendix~\ref{app:core-proofs} contains the full proofs of Theorems~\ref{thmP01}, \ref{thmP1}, \ref{them29}, and~\ref{them215}; the remaining proofs are in Supplementary Appendix~\ref{supp-supp:sec:proofs}.

\subsection{Consistency}
\label{Consis}
We establish estimation consistency and recovery of the population threshold set $\mathcal N_{b_n}$ for the ordinary and debiased estimators. The thresholded error bounds account for the discarded true coefficients, and the analysis also gives consistency of the residual-based noise-variance estimators. The assumptions below specify the spectral, signal, and moment conditions needed for these results.

The assumptions follow the form used in debiased and thresholded Ridge analysis \citep{Zhang2020RidgeRR}, with generator and qualification bounds specifying the regularization class.

Let $\mathbf{P}_n$ denote the orthogonal projector onto
$\operatorname{range}(\mathbf{X}_n^{\mathrm T})$, and write
$\pmb{\beta}_\perp=(\mathbf{I}_p-\mathbf{P}_n)\pmb{\beta}$ for the component of
$\pmb{\beta}$ in $\ker(\mathbf{X}_n)$. With the thin SVD in
Appendix~\ref{app:spectral}, equation~\eqref{SVD},
$\mathbf{P}_n=\mathbf{V}\mathbf{V}^{\mathrm T}$.

Write $b_n=C_b n^{-v_b}$ with constants $C_b,v_b>0$.

\noindent\textbf{Assumptions:}
\begin{enumerate}
    \item[1(a).] Polynomial-growth conditions on the eigenvalues of $\mathbf{X}_n^\mathrm{T}\mathbf{X}_n$: there exist constants $c_\lambda, C_\lambda>0,0<\eta \leq 1 / 2$, such that the positive eigenvalues of $\mathbf{X}_n^\mathrm{T}\mathbf{X}_n$ satisfy the inequality: 
$C_\lambda n \geq \lambda_1 \geq \lambda_2 \geq \cdots \geq \lambda_s \geq c_\lambda n^{2\eta}.$
    \item[1(b).] Euclidean energy of ground truth: $\|\pmb{\beta}\|_2 =\sqrt{\sum_{i=1}^p \beta_i^2} = O\left(n^{\alpha_\beta}\right)$ with $0<\alpha_\beta<(2d-1) \eta$. Here $d$ is the qualification order given in Definition \ref{indexd}.
    \item[2.] The a priori choice of regression parameter:
$\alpha=O\left(n^{2 \eta-\delta-1}\right)$ with a positive constant $\delta$ such that $\frac{\eta+\alpha_\beta}{d}<\delta<2 \eta$. 
    \item[3.] Error structure: The components of the error vector $\pmb{e}_n=[e_1,\cdots,e_n]^\mathrm{T}$ in model \eqref{hgfc} are assumed to be i.i.d., with $\mathbb{E} e_i= 0$, $\mathbb{E} e^2_i= \sigma^2>0$, and $\mathbb{E}\left|e_i\right|^m<\infty$ for some $m>4$.
    \item[4(a).] The dimension of $\pmb{\beta}$ satisfies $p=O\left(n^{\alpha_p}\right)$ for $\alpha_p \in (0, m \eta)$, where $\eta$ and $m$ are defined in Assumptions 1(a) and 3, respectively.
    \item[4(b).] Threshold conditions: $v_b + \frac{\alpha_p}{m} < \eta$. Assume there exists a constant $0<c_b<1$ such that $\max\limits_{i \notin \mathcal{N}_{b_n}}\left|\beta_i\right| \leq c_b b_n$, and $\min\limits_{i \in \mathcal{N}_{b_n}}\left|\beta_i\right| \geq \frac{b_n}{c_b}$.
    \item[5(a).] Polynomial-growth condition for thresholded ground truth: there exists a positive number $\alpha_\sigma\in (0, \eta]$ such that $\sum_{j \notin \mathcal{N}_{b_n}}\left|\beta_j\right|=O\left(n^{v_b-\alpha_\sigma}\right)$ and $\left|\mathcal{N}_{b_n}\right| = O\left(n^{2(\eta-\alpha_\sigma)}\right)$.
    \item[5(b).] Negligibility conditions:
    \begin{enumerate}
    \item[(i)] The discarded coefficients satisfy
    \begin{equation*}
    \max_{1\leq i\leq p}
    \left\{|\beta_i|\,\mathbf{1}_{\{i\notin\mathcal{N}_{b_n}\}}\right\}
    =o\bigl((n\log n)^{-1/2}\bigr).
    \end{equation*}
    \item[(ii)] The null-space component satisfies
    \begin{equation*}
    \|\pmb{\beta}_\perp\|_\infty
    =o\bigl((n\log n)^{-1/2}\bigr).
    \end{equation*}
    \end{enumerate}

\end{enumerate}

\begin{remark}
\label{remarkThreshold}
Assumption 4(b) separates the retained coefficients from those below the
threshold. Assumption 5(a) controls the size of the retained set and the
total magnitude of the discarded coefficients. Assumption 5(b)(i) makes
the discarded coefficients negligible after normalization for Gaussian
approximation, whereas Assumption 5(b)(ii) controls the component of the
true parameter that cannot be identified from $\mathbf{X}_n\pmb{\beta}$.
The latter condition holds automatically when $\mathbf{X}_n$ has full
column rank; for a rank-deficient design, it is an additional restriction
on the true parameter.
\end{remark}

We first establish consistency in the maximum norm
$\|\pmb{\beta}\|_\infty := \max_{1\leq i\leq p}|\beta_i|$ for the
estimators in \eqref{generalLR} and \eqref{debiasLR}. We then study variable
selection and the Euclidean error of their thresholded versions in
\eqref{threshold-beta}, retaining the contribution of the discarded
coefficients explicitly.

\begin{theorem}
\label{thmP01}
Suppose that Assumptions 1--3, 4(a), and 5(b)(ii) hold. Then
$\left\| \hat{\pmb{\beta}}_\alpha - \pmb{\beta} \right\|_\infty = O_p\left( n^{\frac{\alpha_p}{m}-\eta}\right)$.
Moreover,
$\left\| \tilde{\pmb{\beta}}_\alpha - \pmb{\beta} \right\|_\infty  = O_p\left( n^{\frac{\alpha_p}{m}-\eta}\right)$.

\end{theorem}

\begin{proof}[Proof sketch of Theorem~{\upshape\ref{thmP01}}]
The spectral decompositions in Appendix~\ref{app:spectral} separate each
estimation error into a centered noise term, a spectral bias, and the
unidentifiable component $-\pmb\beta_\perp$. Under Assumptions 1--2,
the ordinary and debiased spectral biases have maximum norms
$O(n^{\alpha_\beta-d\delta})$ and $O(n^{\alpha_\beta-2d\delta})$,
respectively. The moment bound in Supplementary Lemma~\ref{supp-lemma1}
gives $O_p(n^{\alpha_p/m-\eta})$ for both noise terms in maximum norm.
Assumptions 2 and 5(b)(ii) make the spectral and null-space biases
negligible on this scale. Appendix~\ref{app:consistency} gives the full proof.
\end{proof}

Theorem~\ref{thmP01} establishes maximum-norm consistency. The threshold
separation in Assumption 4(b) additionally yields recovery of the
population threshold set.

\begin{theorem}
\label{thmP1}
Suppose that Assumptions 1--4 and 5(b)(ii) hold. Then the threshold sets
associated with $\hat{\pmb{\beta}}_\alpha$ and
$\tilde{\pmb{\beta}}_\alpha$ satisfy
\begin{equation}
\label{Prob1}
\begin{aligned}
\mathbb{P}\left(\widehat{\mathcal{N}}_{b_n}\neq\mathcal{N}_{b_n}\right)
&=O\left(n^{-(m\eta-\alpha_p-mv_b)}\right),\\
\mathbb{P}\left(\widetilde{\mathcal{N}}_{b_n}\neq\mathcal{N}_{b_n}\right)
&=O\left(n^{-(m\eta-\alpha_p-mv_b)}\right).
\end{aligned}
\end{equation} 
\end{theorem}

The next result bounds the Euclidean error of the thresholded estimators,
including the contribution of coefficients outside $\mathcal{N}_{b_n}$.

\begin{theorem}\label{Thm1}
Suppose that Assumptions 1--4, 5(a), and 5(b)(ii) hold. Write
$N=\mathcal N_{b_n}$, and let $\pmb\beta_{N^c}$ denote the vector of
discarded coefficients. Then
\begin{align*}
\|\hat{\pmb{\theta}}-\pmb{\beta}\|_2
&=O_p\!\left(n^{-\alpha_\sigma}+\|\pmb\beta_{N^c}\|_2\right),\\
\|\tilde{\pmb{\theta}}-\pmb{\beta}\|_2
&=O_p\!\left(n^{-\alpha_\sigma}+\|\pmb\beta_{N^c}\|_2\right).
\end{align*}
Assumptions 4(b) and 5(a) give
$\|\pmb\beta_{N^c}\|_2=O(n^{-\alpha_\sigma/2})$, so both errors are
$O_p(n^{-\alpha_\sigma/2})$ and converge to zero in probability.
If $\pmb\beta_{N^c}=0$ for all sufficiently large $n$, both rates improve
to $O_p(n^{-\alpha_\sigma})$.
\end{theorem}

Furthermore, we demonstrate the consistency and the rate of convergence of the estimators for $\sigma^2$.

\begin{theorem}
\label{them27}
Suppose that Assumptions 1--4, 5(a), and 5(b)(ii) hold. Then
\begin{equation}
\label{thm27}
\begin{aligned}
\left|\widehat{\sigma}^2-\sigma^2\right|&=O_p\left(n^{-\alpha_\sigma}\right),\\
\left|\widetilde{\sigma}^2-\sigma^2\right|&=O_p\left(n^{-\alpha_\sigma}\right).
\end{aligned}
\end{equation}
\end{theorem}

\subsection{Gaussian approximation and bootstrap inference}
\label{Gauss}
Distributional approximations complement consistency results when constructing confidence regions \citep{van2000asymptotic}. We study Gaussian approximations to the distributions of maximum normalized estimation errors for the regularized and thresholded estimators, and establish bootstrap validity for a priori tuning. The approximating distributions may depend on the sample size. For this purpose, we introduce one additional assumption.

For $i=1,\ldots,p$, define $\tau_i$ by
\begin{equation}
\label{tau}
\tau_i=\sqrt{\sum_{k=1}^s v_{ik}^2[1+r_\alpha(\lambda_k/n)]^2
g_\alpha^2(\lambda_k/n)\frac{\lambda_k}{n^2}+\frac1n}.
\end{equation}

For $1\leq i\leq p$ and $1\leq l\leq n$, write
\begin{equation*}
a_{il}=\frac{1}{n\tau_i}\sum_{k=1}^s v_{ik}u_{lk}g_\alpha(\lambda_k/n)\sqrt{\lambda_k}[1+r_\alpha(\lambda_k/n)].
\end{equation*}
\noindent\textbf{Assumption 6.} One of the following conditions holds:
\begin{enumerate}
\item[(A)] $\displaystyle
\max_{\substack{1\leq i\leq p\\1\leq l\leq n}}|a_{il}|
=o\!\left(\frac{\min\{n^{(\alpha_\sigma-1)/2},\,n^{-1/3}\}}
{(\log n)^{3/2}}\right).
$
\item[(B)] $\alpha_\sigma<1/2$, $p=o(n^{\alpha_\sigma}\log^{-3}(n))$, and
\begin{equation*}
\max_{\substack{1\leq i\leq p\\1\leq l\leq n}}|a_{il}|
=O\!\left(\frac{n^{-\alpha_\sigma}}{(\log n)^{3/2}}\right).
\end{equation*}
\end{enumerate}
According to the error decomposition in \eqref{debiasest}, the quantity
$\sum_{l=1}^n a_{il}e_l$ asymptotically approximates the normalized estimation
error $((\widetilde{\pmb{\beta}}_\alpha)_i-\beta_i)/\tau_i$. Assumption 6 controls
the standardized noise coefficients and, in case (B), the growth of the dimension.

\begin{remark}
Let $\mathbf Q_n$ be an orthogonal matrix with
$\max_{i,l}|(Q_n)_{li}|=O(n^{-1/2})$. For $k\in\{1,2\}$, set
$p=kn$ and $\mathbf X_n=\sqrt{n/k}[\mathbf Q_n,\ldots,\mathbf Q_n]$
with $k$ blocks. All positive eigenvalues are $n$. Under Assumptions
1--2, $w_n=1-r_\alpha(1)^2\to1$, and
\begin{equation*}
a_{i+(j-1)n,l}=\frac{w_n(Q_n)_{li}}{\sqrt{k+w_n^2}},
\quad 1\leq i,l\leq n,\quad 1\leq j\leq k.
\end{equation*}
Hence $\max_{i,l}|a_{il}|=O(n^{-1/2})$ and
$\sum_l a_{il}^2=w_n^2/(k+w_n^2)\to1/(k+1)$, verifying
Assumption 6(A) and \eqref{eq:gaussian-nondegeneracy}.
The case $k=1$ has full column rank, whereas $k=2$ has $s=n<p$;
in the latter case, $\pmb\beta=(\mathbf b^{\mathrm T},\mathbf b^{\mathrm T})^{\mathrm T}$
satisfies $\pmb\beta_\perp=0$.
\end{remark}

With $G_\alpha$ as in \eqref{method-filters}, define for $x\in\mathbf{R}$

\begin{equation*}
H(x)=\mathbb{P}\!\left(\max_i
\frac{\left|\sum\limits_{k=1}^s v_{ik}G_\alpha(\lambda_k/n)
\sqrt{\lambda_k}\xi_k\right|}{n\tau_i}\leq x\right),
\end{equation*}

where $\max_i$ denotes $\max_{1\leq i\leq p}$, and $\xi_k, k=1,2, \cdots, s$ are independent normal random variables with mean 0 and variance $\sigma^2=\mathbb{E} e_1^2$. 
The distribution of the maximum normalized error is approximated by the sample-size-dependent distribution $H$ of a maximum of absolute Gaussian variables.
The following theorem gives the Gaussian approximation for the debiased estimator.

\begin{theorem}\label{them29}
Suppose that Assumptions 1--3, 4(a), 5(a), 5(b)(ii), and 6 hold. Assume that
there exists a constant $c_a>0$, independent of $n$, such that
\begin{equation}
\label{eq:gaussian-nondegeneracy}
\min_{1\leq i\leq p}\sum_{l=1}^n a_{il}^2\geq c_a
\end{equation}
for all sufficiently large $n$. Then, as $n\to\infty$,
\begin{equation*}
\sup_{x\geq0}\left|\mathbb{P}\left(
\max_{1\leq i\leq p}\frac{|(\tilde{\pmb{\beta}}_\alpha)_i-\beta_i|}{\tau_i}
\leq x\right)-H(x)\right|\to0.
\end{equation*} 
\end{theorem}
\begin{proof}[Proof sketch of Theorem~{\upshape\ref{them29}}]
Let $T_n=\max_i|(\tilde{\pmb\beta}_\alpha)_i-\beta_i|/\tau_i$
and $M_n=\max_i|\sum_l a_{il}e_l|$. The spectral decomposition,
weighted Cauchy--Schwarz, and $\tau_i\geq n^{-1/2}$ give
\begin{equation*}
|T_n-M_n|\leq\Delta_n
:=Cn^{\alpha_\beta+\eta-2d\delta}
+\sqrt n\|\pmb\beta_\perp\|_\infty.
\end{equation*}
Assumptions 2 and 5(b)(ii) imply
$\Delta_n\sqrt{\log n}\to0$. Condition~\eqref{eq:gaussian-nondegeneracy}
keeps the coordinate variances bounded away from zero.
Assumption 6 and Supplementary Lemma~\ref{supp-lemma:rectangle-clt}
give a Gaussian approximation even when the covariance is singular;
its boundary bound absorbs the bias. Appendix~\ref{app:gaussian}
gives the full proof.
\end{proof}

Let $c_{1-\alpha^*}$ denote the $(1-\alpha^*)$-quantile of $H$.
Under the assumptions of Theorem~\ref{them29}, the confidence region
\begin{equation*}
\left\{\pmb\beta\in\mathbf R^p:
\max_{1\leq i\leq p}
\frac{|(\tilde{\pmb\beta}_\alpha)_i-\beta_i|}{\tau_i}
\leq c_{1-\alpha^*}\right\}
\end{equation*}
has coverage probability $1-\alpha^*+o(1)$, uniformly for
$\alpha^*\in[\alpha_0^*,\alpha_1^*]$ whenever
$0<\alpha_0^*<\alpha_1^*<1$ are fixed.

Analogously, for $1\leq i\leq p$ and $x\in\mathbf{R}$, define

\begin{align}
\tau_i^*&=\sqrt{\frac{1}{n^2}\sum_{k=1}^s v_{ik}^2
g_\alpha^2(\lambda_k/n)\lambda_k+\frac1n},\label{tau*}\\
H^*(x)&=\mathbb{P}\!\left(\max_i
\frac{\left|\sum\limits_{k=1}^s v_{ik}g_\alpha(\lambda_k/n)
\sqrt{\lambda_k}\xi_k\right|}{n\tau_i^*}\leq x\right).\notag
\end{align}
The following theorem gives the Gaussian approximation for the ordinary estimator.

\begin{theorem}
\label{them210}
Suppose that the assumptions of Theorem~\ref{them29} hold. Then, as $n\to\infty$,
\begin{equation*}
\sup_{x\geq0}\left|\mathbb{P}\left(
\max_{1\leq i\leq p}\frac{|(\hat{\pmb{\beta}}_\alpha)_i-\beta_i|}{\tau^*_i}
\leq x\right)-H^*(x)\right|\to0.
\end{equation*}

\end{theorem}

We next consider Gaussian approximation for the thresholded estimators $\tilde{\pmb{\theta}}$ and $\hat{\pmb{\theta}}$. Assumption 5(b)(i) controls the contribution of the discarded coefficients. Write $N=\mathcal{N}_{b_n}$. For $N\ne\varnothing$ and $x\in\mathbf{R}$, define
\begin{align*}
H_N(x)&=\mathbb{P}\!\biggl(\max_{i\in N}\frac{|\sum_{k=1}^s v_{ik}G_\alpha(\tfrac{\lambda_k}{n})\sqrt{\lambda_k}\xi_k|}{n\tau_i}\leq x\biggr),\\
H_N^*(x)&=\mathbb{P}\!\biggl(\max_{i\in N}\frac{|\sum_{k=1}^s v_{ik}g_\alpha(\tfrac{\lambda_k}{n})\sqrt{\lambda_k}\xi_k|}{n\tau_i^*}\leq x\biggr).
\end{align*}

\begin{theorem}
\label{them215}
Suppose that Assumptions 1--4, 5(a), 5(b), and 6 hold, together with
\eqref{eq:gaussian-nondegeneracy}. Assume that $N\ne\varnothing$ for all
sufficiently large $n$. Then, as $n\to\infty$,
\begin{equation*}
\sup_{x\geq0}\left|\mathbb{P}\left(
\max_{1\leq i\leq p}\frac{|\tilde{\theta}_i-\beta_i|}{\tau_i}
\leq x\right)-H_N(x)\right|\to0,
\end{equation*}
and
\begin{equation*}
\sup_{x\geq0}\left|\mathbb{P}\left(
\max_{1\leq i\leq p}\frac{|\hat{\theta}_i-\beta_i|}{\tau^*_i}
\leq x\right)-H_N^*(x)\right|\to0.
\end{equation*}
\end{theorem}

\begin{proof}[Proof sketch of Theorem~{\upshape\ref{them215}}]
By Theorem~\ref{thmP1}, both estimated threshold sets equal $N$ with
probability tending to one. On this event, the retained coordinates
have the unthresholded estimation errors, while Assumption 5(b)(i)
makes the maximum normalized error over discarded coordinates
$o((\log n)^{-1/2})$. The Gaussian approximation arguments for
Theorems~\ref{them29} and~\ref{them210} apply to the rows indexed by $N$;
Gaussian anti-concentration controls the additional perturbation.
The full proof is given in Appendix~\ref{app:gaussian}.
\end{proof}

The following result establishes the validity
of the unnormalized bootstrap statistic in Algorithm~\ref{algorithm.bootstrap}.
It extends
the bootstrap argument of \citet[Theorem~3]{Zhang2020RidgeRR} to the
ordinary and debiased estimators considered here.

\begin{theorem}[Bootstrap validity]
\label{thm:bootstrap-validity}
Suppose the assumptions of Theorem~\ref{them215} hold and $s=p$. Let
$\alpha=\alpha_n$ and $b_n$ be deterministic sequences held fixed
throughout the bootstrap replicates. For Algorithm~\ref{algorithm.bootstrap},
write
\begin{equation*}
T_n=\|\hat{\pmb\theta}-\pmb\beta\|_\infty,
\qquad
T_n^*=\|\hat{\pmb\theta}^*-\hat{\pmb\theta}\|_\infty.
\end{equation*}
Then
\begin{equation}
\label{eq:bootstrap-consistency}
\sup_{x\geq0}\left|\mathbb P^*(T_n^*\leq x)
-\mathbb P(T_n\leq x)\right|=o_p(1),
\end{equation}
where $\mathbb P^*(\cdot)=\mathbb P(\cdot\mid\mathbf Y_n)$.
For any fixed $0<\alpha^*<1$, let
$\hat c_{1-\alpha^*}^*$ be the conditional $1-\alpha^*$ quantile
of $T_n^*$. Then
\begin{equation}
\label{eq:bootstrap-coverage}
\mathbb P\!\left(T_n\leq\hat c_{1-\alpha^*}^*\right)
\longrightarrow1-\alpha^*.
\end{equation}
The same limit holds with the empirical quantile
$\hat C_{1-\alpha^*}^*$ from $B=B_n\to\infty$ conditionally independent bootstrap
replicates. All conclusions also hold for the debiased branch, with
hats replaced by tildes.
\end{theorem}

\begin{proof}[Proof sketch]
A Gaussian approximation over rectangles handles the unequal coordinate
scales in the unnormalized maximum. The bias and threshold separation
bounds control the bootstrap bias and ensure conditional recovery of
the retained set. Together with variance consistency, they yield a
common Gaussian approximation for the original and bootstrap maxima.
A quantile argument gives \eqref{eq:bootstrap-coverage}.
Supplementary Appendix~\ref{supp-supp:sec:bootstrap} gives the full proof.
\end{proof}

\subsection{Best-worst-case error}
\label{BestWorst}

We complement the statistical results with a deterministic perturbation analysis for a fixed regularization family. Following the rate framework of \citet{boct2021convergence,WANG2024101826}, we relate noise-free bias to the error obtained by optimizing the regression parameter for each perturbed data vector. The two-sided rate relation below requires additional regularity and spectral gain conditions.

Let $\check{\pmb{\beta}}_\alpha(\check{\mathbf{Y}}_n)$ be the estimator \eqref{generalLR} with data $\check{\mathbf{Y}}_n$, and set $\pmb{\beta}_\alpha=\check{\pmb{\beta}}_\alpha(\mathbf{X}_n\pmb{\beta})$.
For $\sigma>0$, let $\bar B_\sigma(\mathbf{X}_n\pmb{\beta})$ denote the closed Euclidean ball in $\mathbf{R}^n$ centered at $\mathbf{X}_n\pmb{\beta}$ with radius $\sqrt n\,\sigma$.
We study the \emph{best-worst-case error}
\begin{equation*}
\sup_{\check{\mathbf{Y}}_n\in\bar B_\sigma(\mathbf{X}_n\pmb{\beta})}
\inf_{\alpha>0}\|\check{\pmb{\beta}}_\alpha(\check{\mathbf{Y}}_n)-\pmb{\beta}\|_2.
\end{equation*}

\begin{definition}
\label{DefRegular}
A generator $g_\alpha$ is called regular if the following conditions, adapted from \citet[Definition~1]{boct2021convergence}, hold:
  \begin{enumerate}[(D3-3)]
    \item[(D3-1)] $\lim_{\alpha\to\infty}|r_\alpha(\lambda)|=1$ for each $\lambda\in(0,+\infty]$, and, for every $\alpha>0$, the function $\lambda\mapsto r_\alpha(\lambda)$ is nonnegative and nonincreasing on $(0,\alpha)$.
    \item[(D3-2)] There exists a monotonically decreasing, continuous function $R_\alpha:(0, \infty) \to[0,1]$ for every $\alpha>0$ such that $R_\alpha \ge\left|r_\alpha\right|$ and $\alpha \mapsto R_\alpha(\lambda)$ is continuous and monotonically increasing for every fixed $\lambda>0$.
    \item[(D3-3)] There exists a constant $c_1\in(0,1)$ such that $R_\alpha(\alpha)\leq c_1$ for all $\alpha>0$.
  \end{enumerate}
\end{definition}

\begin{lemma}
\label{BestWorst2}
Let $g_\alpha$ satisfy Definition~\ref{DefRegular} for all $\alpha>0$, uniformly in $n$. Suppose $\pmb{\beta}\ne0$, $r_\alpha\ge0$ on the positive spectrum $\operatorname{Spec}^{+}(\mathbf A_n)$ for all $\alpha>0$, and $\alpha\mapsto\|\pmb{\beta}_\alpha-\pmb{\beta}\|_2$ is continuous and nondecreasing. Then each $\sigma>0$ determines a unique $\alpha(\sigma)>0$ satisfying $\sqrt{\alpha(\sigma)}\|\pmb{\beta}_{\alpha(\sigma)}-\pmb{\beta}\|_2=\sigma$, and
\begin{equation}
\label{eqHXY18}
\sup_{\check{\mathbf Y}_n\in\bar B_\sigma(\mathbf X_n\pmb{\beta})}
\inf_{\alpha>0}
\|\check{\pmb{\beta}}_\alpha(\check{\mathbf Y}_n)-\pmb{\beta}\|_2
\le \frac{C_1\sigma}{\sqrt{\alpha(\sigma)}}.
\end{equation}
For fixed $\kappa>0$, if
\begin{equation}
\label{bw-gain}
\max_{\lambda\in\operatorname{Spec}^{+}(\mathbf A_n)}
\inf_{0<\alpha\le\alpha(\sigma)}
\sqrt{\alpha(\sigma)\lambda}\,g_\alpha(\lambda)\ge\kappa,
\end{equation}
then
\begin{equation}
\label{eqHXY19}
\sup_{\check{\mathbf Y}_n\in\bar B_\sigma(\mathbf X_n\pmb{\beta})}
\inf_{\alpha>0}
\|\check{\pmb{\beta}}_\alpha(\check{\mathbf Y}_n)-\pmb{\beta}\|_2
\ge \frac{C_2\sigma}{\sqrt{\alpha(\sigma)}}.
\end{equation}
Here $C_1,C_2>0$ are independent of $n$ and $\sigma$.
\end{lemma}

Using the rate transformation in \citet[Definition~5 and Lemma~6]{boct2021convergence} and \citet{WANG2024101826}, Lemma~\ref{BestWorst2} yields the following equivalence on the scales satisfying \eqref{bw-gain}.

\begin{theorem}
\label{thmBW}
Under the assumptions of Lemma~\ref{BestWorst2}, let $\varphi:[0,\infty)\to[0,\infty)$ be continuous, strictly increasing and $\zeta$-homogeneous\footnote{A function $\varphi$ is $\zeta$-homogeneous if there exists an increasing function $\zeta:[0,\infty)\to[0,\infty)$ such that $\varphi(\gamma\alpha)\le\zeta(\gamma)\varphi(\alpha)$ for all $\alpha,\gamma>0$.}, with $\varphi(0)=0$. Set $\check\varphi(\alpha)=\sqrt\alpha\,\varphi(\alpha)$ and $\psi(\sigma)=\sigma/\sqrt{\check\varphi^{-1}(\sigma)}$. For fixed $\kappa>0$, independent of $n$, the following statements are equivalent:
\begin{enumerate}
\item\label{bw-noisy-rate} There exists $c>0$ such that, for all $(n,\sigma)$ satisfying \eqref{bw-gain},
\begin{equation*}
\sup_{\check{\mathbf Y}_n\in\bar B_\sigma(\mathbf X_n\pmb{\beta})}
\inf_{\alpha>0}
\|\check{\pmb{\beta}}_\alpha(\check{\mathbf Y}_n)-\pmb{\beta}\|_2
\le c\,\psi(\sigma).
\end{equation*}
\item\label{bw-bias-rate} There exists $\check c>0$ such that $\|\pmb{\beta}_\alpha-\pmb{\beta}\|_2\le\check c\,\varphi(\alpha)$ for all $(n,\alpha)$ satisfying
\begin{equation}
\label{bw-gain-alpha}
\max_{\lambda\in\operatorname{Spec}^{+}(\mathbf A_n)}
\inf_{0<t\le\alpha}
\sqrt{\alpha\lambda}\,g_t(\lambda)\ge\kappa.
\end{equation}
\end{enumerate}
The constants $c,\check c$ are independent of $n$, $\sigma$ and $\alpha$.
\end{theorem}

If the bias bound holds for all $(n,\alpha)$, the noisy bound holds for all $(n,\sigma)$ without \eqref{bw-gain}.

\begin{remark}
Lemma~\ref{BestWorst2} and Theorem~\ref{thmBW} also apply to the debiased
estimator when its generator $G_\alpha=(1+r_\alpha)g_\alpha$ satisfies
their assumptions, with $\alpha(\sigma)$ defined using the corresponding
noise-free estimator.
\end{remark}

\section{Numerical experiments}
\label{simulation}
We examine sparse and non-sparse regression and an inverse source problem. The implementation and tuning rules are specified below; the computational environment is described in Supplementary Section~\ref{supp-supp:sec:computational-environment}.

The regression experiments compare ordinary, debiased, and thresholded estimators under the tuning rules specified below. Details of the experimental scope and regression-bootstrap and independent-validation protocols are given in Supplementary Section~\ref{supp-supp:sec:regression-settings}. The inverse source experiment uses the observation and resampling scheme described in Section~\ref{sec:isp}.

\subsection{Sparse case}\label{Simu1}
We generate the design matrix $\mathbf{X}_n$, coefficient vector $\pmb{\beta}$, and error vector $\pmb{e}_n$ as follows.
\begin{itemize}
    \item Design matrix $\mathbf{X}_n$: Define $\mathbf{X}_n=\left[x_1, \cdots, x_n\right]^\mathrm{T}$ with $x_i=\left(x_{i 1}, \cdots, x_{i p}\right)^\mathrm{T} \in \mathbf{R}^p, i=1, \cdots, n$. Generate $x_1, x_2, \cdots$ as i.i.d. normal random vectors with mean $\pmb{0}$ and covariance matrix $\Sigma \in$ $\mathbf{R}^{p \times p}$. We select $\Sigma$ with diagonal elements equal to 2.0 and off-diagonal elements equal to 0.5. 
    Case I uses a fixed Gaussian realization generated as described above. Normal and Laplace samples share the design and coefficients within each case.
    \item Coefficient vector $\pmb{\beta}$: Generate a zero vector $\pmb{\beta}=\left(\beta_1, \cdots, \beta_p\right)^\mathrm{T} (p\geq 20)$, then randomly select 20 components of $\pmb{\beta}$. Assign values as follows: set five components to 2, five to -2, five to 1, and five to -1.
    \item Error vector $\pmb{e}_n$: Use centered normal errors with variance 4 or centered Laplace errors with scale $\sqrt{2}$, also giving variance 4.
\end{itemize}

The response is generated according to model~\eqref{hgfc}. For Lasso and Ridge, the regression parameter minimizes the ordinary, unthresholded coefficient error $\|\hat{\pmb{\beta}}_{\alpha}-\pmb{\beta}\|$ over the reported grids. SC similarly selects its truncation rank using the true coefficients. Thresholds are selected subsequently for each available estimator branch; the regression parameter and threshold are not jointly optimized.

For the seven iterative regression methods (Landweber, Showalter, SOAR, HBF, FAR, AR$^\kappa$, and Nesterov), we use the truncated discrepancy principle (TDP). The output is $\hat{\pmb{\beta}}=\pmb{\beta}_{k_0}$, where $k_0=\min(k_{\max},k^*)$, and $k^*$ is the first index satisfying the discrepancy condition. The right-hand inequality below holds for every integer $k$ with $1\leq k<k^*$:
    \begin{equation}
\label{discrepancy}
\|\mathbf{Y}_n-\mathbf{X}_n\pmb{\beta}_{k^*}\|
\leq\varsigma\|\pmb{e}_n\|<\|\mathbf{Y}_n-\mathbf{X}_n\pmb{\beta}_k\|,
\end{equation}
Here $\varsigma>0$ is specified for each experiment.

\begin{remark}
The normalized estimator in \eqref{generalLR} can be implemented by supplying $\bar{\mathbf{X}}_n=\mathbf{X}_n/\sqrt n$ and $\bar{\mathbf{Y}}_n=\mathbf{Y}_n/\sqrt n$ to the corresponding regression algorithm, since
$g_\alpha(\bar{\mathbf{X}}_n^\mathrm{T}\bar{\mathbf{X}}_n)\bar{\mathbf{X}}_n^\mathrm{T}\bar{\mathbf{Y}}_n=n^{-1}g_\alpha(\mathbf{X}_n^\mathrm{T}\mathbf{X}_n/n)\mathbf{X}_n^\mathrm{T}\mathbf{Y}_n$.
The experiments use unscaled inputs. Relating their numerical parameters to the normalized formulation requires a method-specific conversion; for Landweber, $h=n\Delta t$, as explained after \eqref{debiasedLandweber}. Step sizes, coefficients, iteration limits, and search grids are given in Supplementary Section~\ref{supp-supp:sec:regression-settings} (Online Resource~1).
\end{remark}
The two cases have different $p/n$ ratios, with normal and Laplace noise realizations paired across the designs. The sparse point-estimation and regression-bootstrap experiments use the common numerical settings in Supplementary Table~\ref{supp-supp:tab:regression-parameters}.

\subsubsection{Case I: \texorpdfstring{$n=p=1000$}{n=p=1000}}
For both sparse cases, the TDP multiplier is $\varsigma=1.01$ for SOAR and $\varsigma=0.9$ for the other six iterative methods. SOAR uses $s^*=1.5$. All seven iterative methods reach the discrepancy condition before their respective caps in these four samples.

For both cases, coordinate thresholds minimize the $L_2$ error separately for the ordinary and debiased estimates using the known coefficient vector. Each candidate threshold is applied to an unchanged fitted vector. The grid and tie-breaking rule are specified in Supplementary Section~\ref{supp-supp:sec:regression-settings}. Here $k_r$ denotes the number of singular directions retained by SC.

\begin{table*}[t]
\centering
\caption{Single-sample $L_2$ estimation errors and selected thresholds for Case I ($n=p=1000$), using the fixed Gaussian design.}
\label{case3}
\begingroup
\fontsize{8}{8.7}\selectfont
\setlength{\tabcolsep}{3.5pt}
\renewcommand{\arraystretch}{0.93}
\begin{tabular}{@{}lc S[table-format=3.4] S[table-format=1.4] S[table-format=3.4] *{3}{S[table-format=1.4]} @{}}
\toprule
\multirow{2}{*}{Method} & \multirow{2}{*}{Parameter} & \multicolumn{2}{c}{Unthresholded $L_2$ error} & \multicolumn{2}{c}{Thresholded $L_2$ error} & \multicolumn{2}{c}{Threshold $b_n$} \\
\cmidrule(lr){3-4}\cmidrule(lr){5-6}\cmidrule(l){7-8}
 & & {Ord.} & {Deb.} & {Ord.} & {Deb.} & {Ord.} & {Deb.} \\
\midrule
\multicolumn{8}{@{}l}{\textit{Normal errors}} \\[-0.2ex]
Landweber & 589 & 3.6752 & 3.3431 & 2.4142 & 1.6994 & 0.3577 & 0.4160 \\
Showalter & 590 & 3.6745 & 3.3427 & 2.4136 & 1.6986 & 0.3577 & 0.4160 \\
HBF & 176 & 3.6872 & 3.3511 & 1.9931 & 1.6642 & 0.3953 & 0.4213 \\
$\textrm{AR}^{\kappa}$ & 263 & 3.6974 & 3.3632 & 2.3125 & 1.7152 & 0.3745 & 0.4177 \\
SOAR & 413 & 3.8433 & 3.4443 & 2.3990 & 1.8841 & 0.4180 & 0.4724 \\
Nesterov & 89 & 3.7102 & 3.3732 & 2.3130 & 1.7293 & 0.3758 & 0.4173 \\
FAR & 152 & 3.6596 & 3.3250 & 1.8095 & 1.6146 & 0.4021 & 0.4204 \\
LS & \textemdash & 292.2430 & {\textemdash} & 292.1820 & {\textemdash} & 0.9995 & {\textemdash} \\
SC & $k_r=865$ & 3.5080 & 3.5080 & 0.9157 & 0.9157 & 0.3590 & 0.3590 \\
Lasso & $\alpha=0.133$ & 0.4302 & {\textemdash} & 0.3785 & {\textemdash} & 0.4600 & {\textemdash} \\
Ridge & $\alpha=91.1$ & 3.1742 & 3.3684 & 1.4581 & 0.8854 & 0.4010 & 0.4484 \\
\midrule
\multicolumn{8}{@{}l}{\textit{Laplace errors}} \\[-0.2ex]
Landweber & 563 & 3.7400 & 3.3898 & 2.6023 & 1.8757 & 0.3834 & 0.4151 \\
Showalter & 563 & 3.7405 & 3.3899 & 2.6039 & 1.8763 & 0.3834 & 0.4148 \\
HBF & 169 & 3.7253 & 3.3987 & 2.1572 & 1.8441 & 0.4290 & 0.4188 \\
$\textrm{AR}^{\kappa}$ & 259 & 3.7600 & 3.4123 & 2.5208 & 1.9002 & 0.3978 & 0.4168 \\
SOAR & 402 & 3.9001 & 3.4981 & 2.6113 & 2.0730 & 0.4266 & 0.4770 \\
Nesterov & 87 & 3.7638 & 3.4186 & 2.5029 & 1.9046 & 0.4008 & 0.4165 \\
FAR & 144 & 3.6973 & 3.3780 & 2.0393 & 1.8067 & 0.4329 & 0.4156 \\
LS & \textemdash & 119.7841 & {\textemdash} & 119.5001 & {\textemdash} & 0.9960 & {\textemdash} \\
SC & $k_r=855$ & 3.6516 & 3.6516 & 1.1496 & 1.1496 & 0.3650 & 0.3650 \\
Lasso & $\alpha=0.120$ & 0.6343 & {\textemdash} & 0.5722 & {\textemdash} & 0.4470 & {\textemdash} \\
Ridge & $\alpha=91.6$ & 3.2264 & 3.4289 & 1.5939 & 1.0171 & 0.3652 & 0.3916 \\
\bottomrule
\end{tabular}

\endgroup
\end{table*}

Tables~\ref{case3} and~\ref{case2} compare eleven methods under normal and Laplace errors. Ord./Deb. denote ordinary/debiased estimators. The parameter column gives $k_0$ for iterative methods, $k_r$ for SC, and $\alpha$ for Lasso/Ridge; a dash denotes an inapplicable quantity.

\begin{figure*}[tp]
\centering
\includegraphics[width=\textwidth]{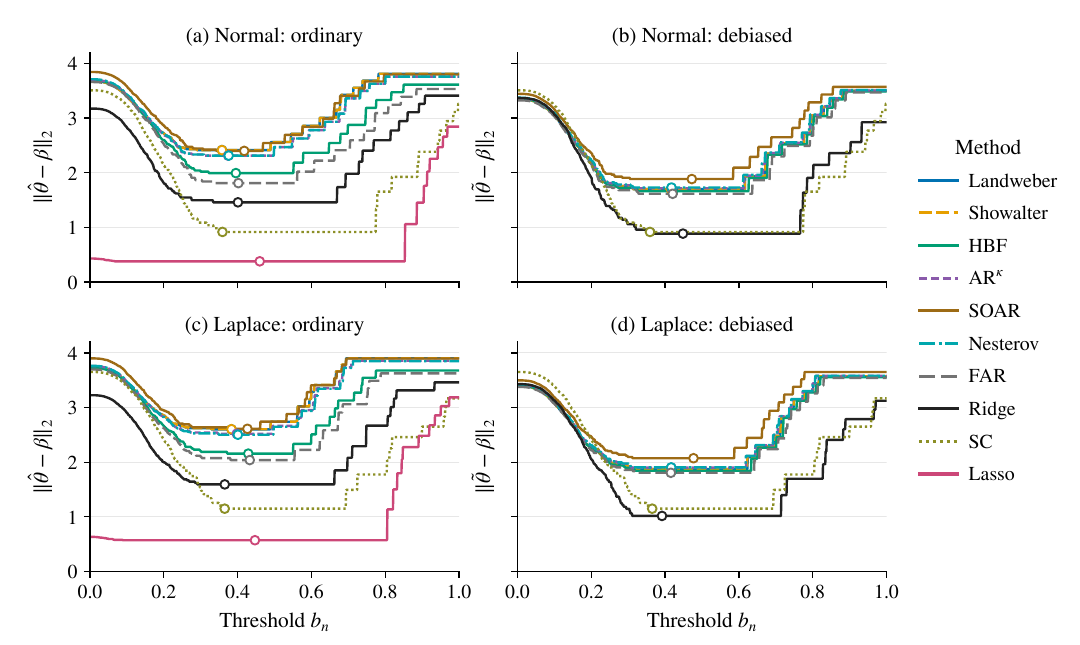}
\caption{Single-sample $L_2$ estimation error versus the threshold in Case I ($n=p=1000$). Panels (a)--(b) use normal errors and (c)--(d) Laplace errors. The left and right columns show the $L_2$ errors of the ordinary and debiased thresholded estimators, respectively. Open circles indicate the selected thresholds. LS is omitted, and Lasso is shown only in the ordinary panels}
\label{fig:case3}
\end{figure*}

In Case I, Lasso has the smallest $L_2$ error at the selected threshold in both samples, followed by debiased Ridge. FAR has the smallest error among the seven iterative methods in both ordinary and debiased branches. LS has much larger errors and is omitted from the plots; its results remain in Table~\ref{case3}.

Figure~\ref{fig:case3} shows the corresponding threshold-error curves. The SC curves coincide for the ordinary and debiased estimators because debiasing leaves its spectral projection unchanged.

Table~\ref{Coverage} compares nine methods in the Case I Gaussian setting. Coverage~I and~II refer to confidence regions centered at the thresholded ordinary and debiased estimators, respectively. LS and Lasso are excluded from this bootstrap comparison.

\begin{table*}[t]
\centering
\caption{Regression-bootstrap calibration for Case I with normal errors ($n=p=1000$, $N=1000$, $B=500$; nominal coverage $95\%$). Error entries are simulation means.}
\label{Coverage}
\begingroup
\fontsize{8}{8.7}\selectfont
\setlength{\tabcolsep}{3pt}
\renewcommand{\arraystretch}{0.98}
\begin{tabular}{@{}l*{9}{S[table-format=2.4]}@{}}
\toprule
Metric & {Landweber} & {Showalter} & {HBF} & {$\mathrm{AR}^{\kappa}$} & {SOAR} & {Nesterov} & {FAR} & {SC} & {Ridge} \\
\midrule
\addlinespace[1.5pt]
Coverage I (\%) & \multicolumn{1}{S[table-format=2.1]}{0.0} & \multicolumn{1}{S[table-format=2.1]}{0.0} & \multicolumn{1}{S[table-format=2.1]}{50.7} & \multicolumn{1}{S[table-format=2.1]}{0.8} & \multicolumn{1}{S[table-format=2.1]}{0.4} & \multicolumn{1}{S[table-format=2.1]}{1.0} & \multicolumn{1}{S[table-format=2.1]}{81.3} & \multicolumn{1}{S[table-format=2.1]}{92.2} & \multicolumn{1}{S[table-format=2.1]}{92.5} \\
Coverage II (\%) & \multicolumn{1}{S[table-format=2.1]}{95.0} & \multicolumn{1}{S[table-format=2.1]}{95.0} & \multicolumn{1}{S[table-format=2.1]}{94.9} & \multicolumn{1}{S[table-format=2.1]}{95.1} & \multicolumn{1}{S[table-format=2.1]}{94.9} & \multicolumn{1}{S[table-format=2.1]}{95.0} & \multicolumn{1}{S[table-format=2.1]}{95.4} & \multicolumn{1}{S[table-format=2.1]}{92.6} & \multicolumn{1}{S[table-format=2.1]}{94.9} \\
\addlinespace[1.5pt]
$L_2$ error I & 2.5006 & 2.5008 & 2.0636 & 2.4098 & 2.5030 & 2.3937 & 1.9143 & 1.0084 & 1.4423 \\
$L_2$ error II & 1.7729 & 1.7726 & 1.7382 & 1.7941 & 1.9692 & 1.7986 & 1.6953 & 1.0084 & 0.8798 \\
\addlinespace[1.5pt]
$L_\infty$ error I & 0.8327 & 0.8328 & 0.7212 & 0.8080 & 0.8243 & 0.8042 & 0.6802 & 0.4364 & 0.5518 \\
$L_\infty$ error II & 0.6553 & 0.6552 & 0.6510 & 0.6620 & 0.7066 & 0.6637 & 0.6397 & 0.4364 & 0.3937 \\
\addlinespace[1.5pt]
Variance error I & 9.2466 & 9.2484 & 6.2483 & 8.5673 & 9.2295 & 8.4523 & 5.3417 & 1.4794 & 3.0519 \\
Variance error II & 4.6071 & 4.6053 & 4.4266 & 4.7185 & 5.6997 & 4.7426 & 4.2099 & 1.4794 & 1.1201 \\
\bottomrule
\end{tabular}
\par\vspace{2pt}
\begin{minipage}{\textwidth}
\fontsize{7.5}{8.5}\selectfont
I/II: ordinary/debiased thresholded estimators. Variance error is the absolute error in estimating $\sigma^2$.
\end{minipage}
\endgroup
\end{table*}

\begin{table*}[t]
\centering
\caption{Single-sample $L_2$ estimation errors and selected thresholds for Case II ($n=1000$, $p=1500$), using the conditional extension of the Case I design.}
\label{case2}
\begingroup
\fontsize{8}{8.7}\selectfont
\setlength{\tabcolsep}{3.5pt}
\renewcommand{\arraystretch}{0.93}
\begin{tabular}{@{}lc *{6}{S[table-format=1.4]} @{}}
\toprule
\multirow{2}{*}{Method} & \multirow{2}{*}{Parameter} & \multicolumn{2}{c}{Unthresholded $L_2$ error} & \multicolumn{2}{c}{Thresholded $L_2$ error} & \multicolumn{2}{c}{Threshold $b_n$} \\
\cmidrule(lr){3-4}\cmidrule(lr){5-6}\cmidrule(l){7-8}
 & & {Ord.} & {Deb.} & {Ord.} & {Deb.} & {Ord.} & {Deb.} \\
\midrule
\multicolumn{8}{@{}l}{\textit{Normal errors}} \\[-0.2ex]
Landweber & 412 & 4.6651 & 4.4778 & 3.4936 & 2.8886 & 0.3104 & 0.3817 \\
Showalter & 413 & 4.6645 & 4.4775 & 3.4929 & 2.8876 & 0.3104 & 0.3817 \\
HBF & 143 & 4.6571 & 4.4676 & 2.9218 & 2.8008 & 0.3485 & 0.4026 \\
$\textrm{AR}^{\kappa}$ & 238 & 4.6762 & 4.4852 & 3.4495 & 2.8993 & 0.3049 & 0.3792 \\
SOAR & 351 & 4.7547 & 4.5191 & 3.4885 & 3.0142 & 0.2888 & 0.3704 \\
Nesterov & 92 & 4.6786 & 4.4882 & 3.3962 & 2.8986 & 0.3140 & 0.3789 \\
FAR & 113 & 4.6592 & 4.4720 & 3.0366 & 2.8265 & 0.3436 & 0.3945 \\
LS & \textemdash & 4.6568 & {\textemdash} & 2.2310 & {\textemdash} & 0.5005 & {\textemdash} \\
SC & $k_r=943$ & 4.5628 & 4.5628 & 2.4123 & 2.4123 & 0.4750 & 0.4750 \\
Lasso & $\alpha=0.141$ & 0.4677 & {\textemdash} & 0.4006 & {\textemdash} & 0.4650 & {\textemdash} \\
Ridge & $\alpha=152.3$ & 4.4379 & 4.5161 & 2.8352 & 2.3766 & 0.3964 & 0.4608 \\
\midrule
\multicolumn{8}{@{}l}{\textit{Laplace errors}} \\[-0.2ex]
Landweber & 391 & 4.7273 & 4.5258 & 3.6797 & 3.0656 & 0.2894 & 0.3415 \\
Showalter & 391 & 4.7278 & 4.5259 & 3.6816 & 3.0664 & 0.2894 & 0.3415 \\
HBF & 138 & 4.6891 & 4.5134 & 3.0928 & 2.9764 & 0.3216 & 0.3490 \\
$\textrm{AR}^{\kappa}$ & 235 & 4.7315 & 4.5313 & 3.6213 & 3.0718 & 0.2934 & 0.3410 \\
SOAR & 341 & 4.8117 & 4.5696 & 3.6743 & 3.1959 & 0.2636 & 0.3438 \\
Nesterov & 89 & 4.7393 & 4.5369 & 3.5999 & 3.0838 & 0.2930 & 0.3393 \\
FAR & 107 & 4.7042 & 4.5205 & 3.2671 & 3.0148 & 0.3202 & 0.3416 \\
LS & \textemdash & 4.7194 & {\textemdash} & 2.4125 & {\textemdash} & 0.4405 & {\textemdash} \\
SC & $k_r=950$ & 4.6262 & 4.6262 & 2.5244 & 2.5244 & 0.3480 & 0.3480 \\
Lasso & $\alpha=0.124$ & 0.6758 & {\textemdash} & 0.5997 & {\textemdash} & 0.4465 & {\textemdash} \\
Ridge & $\alpha=157.5$ & 4.4997 & 4.5767 & 2.9975 & 2.5491 & 0.3460 & 0.3872 \\
\bottomrule
\end{tabular}

\endgroup
\end{table*}

Original-sample Ridge and SC parameters and branch-specific thresholds follow the sparse selection rules above. In Table~\ref{Coverage}, all nine methods share the fixed Case I design, coefficient vector, and outer-noise realizations; outer-sample and bootstrap noises use separate reproducible streams. Case II uses the same numerical settings and threshold-selection rules.

Debiasing reduces the mean thresholded coefficient errors and mean absolute variance-estimation errors for Ridge and all seven iterative methods, whose Coverage~II values range from $94.9\%$ to $95.4\%$. SC has Coverage~II of $92.6\%$. Its ordinary and debiased fits, coefficient errors, and variance estimates are identical; the $0.4$ percentage-point coverage difference results from separately generated finite bootstrap samples. An independent validation with fresh noise realizations is reported in Supplementary Table~\ref{supp-supp:tab:regression-independent-validation}; its coverage Monte Carlo standard errors are given in Supplementary Table~\ref{supp-supp:tab:regression-bootstrap-mcse}.

\subsubsection{Case II: \texorpdfstring{$n=1000,\; p=1500$}{n=1000, p=1500}}
Case II uses the same method settings and threshold-selection rules as Case I, with $n=1000$ and $p=1500$. It retains the first 1000 design columns and coefficients. The additional 500 columns are generated conditionally on the retained design under a Gaussian model whose row covariance has entries $\Sigma_{jj}=2$ and $\Sigma_{jk}=0.5$ for $j\ne k$, $1\leq j,k\leq1500$. The appended coefficients are zero. Supplementary Section~\ref{supp-supp:sec:regression-settings} gives the exact conditional sampling rule and experimental scope. This paired extension preserves the signal and noise and provides an empirical comparison of estimation errors after adding inactive covariates.

Lasso has the smallest error at the selected threshold in both Case II samples, and HBF has the smallest error among the seven iterative methods in both estimator branches. The threshold-error curves are shown in Supplementary Figure~\ref{supp-fig:case2}.

Across the four sparse samples, thresholding lowers the $L_2$ error of every available estimator branch, and debiasing lowers both unthresholded errors and errors at the selected thresholds for all seven iterative methods. Ridge has larger unthresholded errors after debiasing but smaller errors after subsequent thresholding. SC is unchanged by debiasing.


\subsection{Non-sparse case}\label{Simu2}
We compare the regularized estimators $\hat{\pmb{\beta}}_{\alpha}$ and their debiased counterparts $\tilde{\pmb{\beta}}_{\alpha}$ with the same selected regression parameters and numerical settings, without coordinate thresholding.

The experiment uses $n=p=1000$ and the fixed Case I design. The coefficients are drawn independently from the uniform distribution on $[-2,2]$, and the errors are independent $N(0,4)$ variables.
 
The seven iterative methods use TDP in \eqref{discrepancy} with $\varsigma=0.6$. The stopping index selected from the ordinary fit is reused for the debiased fit. All seven methods satisfy the discrepancy condition before their iteration caps.

For Lasso and Ridge, we search the penalty grids in descending order and select the first parameter value for which the residual norm is no greater than $\varsigma\|\pmb e_n\|$. SC selects its truncation rank by minimizing the coefficient error using the true coefficients.
The method-specific numerical settings and penalty grids are given in Supplementary Section~\ref{supp-supp:sec:regression-settings}.

Table~\ref{nonsparse} reports the ordinary and debiased $L_2$ errors and selected iteration counts. Debiasing reduces the errors of Ridge and all seven iterative methods without thresholding. HBF has the smallest ordinary and debiased errors among the iterative methods. HBF and FAR have smaller debiased errors than Ridge; SC is unchanged by debiasing.

Because the sparse and non-sparse comparisons use different penalty-selection rules, their relative performances cannot be attributed to sparsity alone.
\begin{table}[htbp]
\centering
\caption{Single-sample non-sparse coefficient $L_2$ errors and iteration counts ($n=p=1000$). The seven iterative methods use TDP; LS and SC are reference methods. A dash denotes an inapplicable quantity.}
\label{nonsparse}
\begingroup
\fontsize{8}{9}\selectfont
\setlength{\tabcolsep}{4pt}
\renewcommand{\arraystretch}{0.98}
\begin{tabular}{@{}l S[table-format=3.4] S[table-format=4.0] S[table-format=2.4] @{}}
\toprule
Method & {Ordinary $L_2$} & {$k_0$} & {Debiased $L_2$} \\
\midrule
LS & 103.3936 & {\textemdash} & {\textemdash} \\
SC & 7.5771 & {\textemdash} & 7.5771 \\
Lasso & 10.1699 & {\textemdash} & {\textemdash} \\
Ridge & 8.5080 & {\textemdash} & 7.4152 \\
Landweber & 9.4603 & 6296 & 8.3690 \\
Showalter & 9.4604 & 6296 & 8.3690 \\
HBF & 7.8750 & 1390 & 7.2271 \\
$\textrm{AR}^\kappa$ & 9.4503 & 552 & 8.4215 \\
SOAR & 9.4362 & 837 & 8.4205 \\
Nesterov & 9.4570 & 287 & 8.4356 \\
FAR & 9.0519 & 1120 & 7.3342 \\
\bottomrule
\end{tabular}
\endgroup
\end{table}

\subsection{Inverse source problems (ISP) in partial differential equations (PDEs)}
\label{sec:isp}
We illustrate the reconstruction and inference procedures using an elliptic inverse source problem. Given Dirichlet and Neumann boundary data $q_1$ and $q_2$ on $\Gamma=\partial\Omega$, the aim is to recover a source $f$ supported in $\Omega_0\subset\Omega$ from
\begin{equation}
\label{inversesource}
\left\{
\begin{array}{l}
-\Delta u+u=f\chi_{\Omega_0}\quad\text{in }\Omega,\\
u=q_1\quad\text{and}\quad \dfrac{\partial u}{\partial\mathbf n}=q_2
\quad\text{on }\Gamma.
\end{array}
\right.
\end{equation}
Here $\mathbf n$ is the unit outward normal and $\chi_{\Omega_0}$ is the indicator of the source region. Finite element discretization gives a linear dependence of the boundary response on the source coefficients, motivating the use of the regression methods in Section~\ref{example}. We assess their finite-sample reconstruction errors and empirical coverage in this discretized setting.

\begin{figure*}[!t]
\centering
\includegraphics[width=\textwidth]{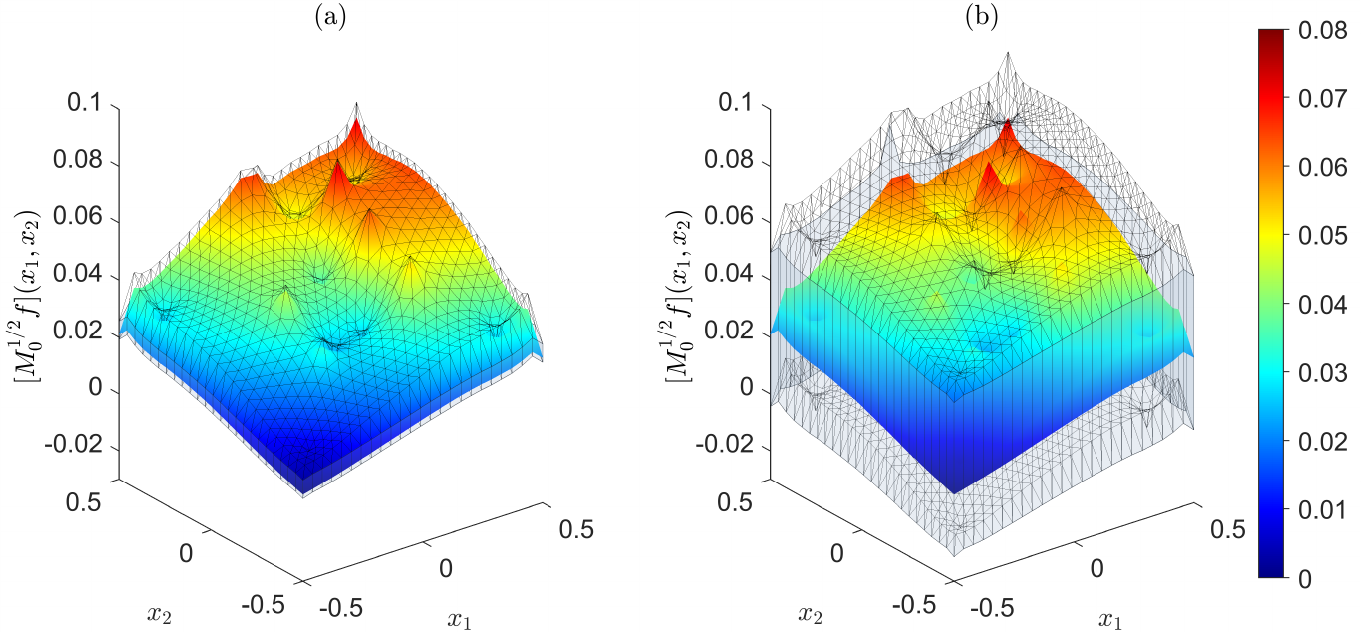}
\caption{Nominal $95\%$ simultaneous confidence regions for the transformed source $M_0^{1/2}f$ using Landweber in the first outer simulation: (a) debiased estimator and (b) ordinary estimator. Here $M_0^{1/2}$ is the full principal square root of the source-region finite element mass matrix. The colored surface interpolates the transformed true-source vector $M_0^{1/2}f^\dagger$, and black wireframes interpolate the lower and upper nodal bounds. Empirical coverage over $N=1000$ simulations is reported in Table~\ref{CoveragePDE}}
\label{fig:pde-landweber}
\end{figure*}

Following \citet{Zhang2018IP}, we take
$\Omega=\{(x_1,x_2):x_1^2+x_2^2<1\}$,
$\Omega_0=(-1/2,1/2)^2$, and
$f^\dagger(x_1,x_2)=(1+x_1+x_2)\chi_{\Omega_0}$.
Exact boundary data are generated with $q_2=0$, and Gaussian noise with standard deviation $\sigma=0.002$ is added to the boundary observations. We compare seven iterative methods using the weighted, effective-noise version of the truncated discrepancy principle in Supplementary equation~\eqref{supp-pde-stopping}. For each method, we use $N=1000$ independent simulations and $B=500$ bootstrap replicates per simulation to construct confidence regions with nominal coverage $95\%$. In each bootstrap sample, the first-stage stopping index is reselected and then held fixed for the second fit used in debiasing.
Supplementary Section~\ref{supp-supp:sec:pde} gives the numerical implementation details (Online Resource~1). Figure~\ref{fig:pde-landweber} shows the single-realization confidence regions for Landweber; corresponding illustrations for the other six methods are provided in Supplementary Figures~\ref{supp-fig:Showalter}--\ref{supp-fig:FAR}.

\begin{table*}[t]
\centering
\caption{Reconstruction and variance-estimation errors and empirical coverage in the inverse source experiment. Coverage I and II refer to $\hat f$ and $\tilde f$, respectively. Overbars denote means over $N=1000$ outer simulations, with $B=500$ bootstrap replicates and nominal coverage $95\%$. The error measures are defined in Supplementary equation~\eqref{supp-pde-error-measures}.}
\label{CoveragePDE}
\begingroup
\footnotesize
\setlength{\tabcolsep}{3pt}
\begin{tabular}{c|c|c|c|c|c|c|c}
\hline
 & Landweber & Showalter & HBF & $\mathrm{AR}^{\kappa}$ & SOAR & Nesterov & FAR \\ \hline
$\varsigma$ & 3.25 & 3.25 & 3.2 & 4.1 & 3.4 & 3.3 & 3.35 \\
Coverage I & 100.0\% & 100.0\% & 100.0\% & 100.0\% & 100.0\% & 100.0\% & 100.0\% \\
Coverage II & 95.4\% & 97.1\% & 94.0\% & 95.0\% & 94.1\% & 94.4\% & 96.2\% \\
$\overline{E_2(\hat f)}$ & 0.100723 & 0.101348 & 0.099999 & 0.088982 & 0.103220 & 0.098893 & 0.100527 \\
$\overline{E_2(\tilde f)}$ & 0.035405 & 0.035606 & 0.035177 & 0.034954 & 0.036332 & 0.034838 & 0.035323 \\
$\overline{E_\infty(\hat f)}$ & 0.009752 & 0.009811 & 0.009683 & 0.008320 & 0.009906 & 0.009585 & 0.009891 \\
$\overline{E_\infty(\tilde f)}$ & 0.003366 & 0.003396 & 0.003328 & 0.003577 & 0.003485 & 0.003287 & 0.003352 \\
$\overline{\lvert\widehat{\sigma}^2-\sigma^2\rvert}$ & $3.725\times10^{-5}$ & $3.777\times10^{-5}$ & $3.667\times10^{-5}$ & $3.994\times10^{-5}$ & $4.102\times10^{-5}$ & $3.587\times10^{-5}$ & $3.969\times10^{-5}$ \\
$\overline{\lvert\widetilde{\sigma}^2-\sigma^2\rvert}$ & $2.068\times10^{-6}$ & $2.131\times10^{-6}$ & $1.994\times10^{-6}$ & $2.474\times10^{-6}$ & $2.352\times10^{-6}$ & $1.926\times10^{-6}$ & $2.051\times10^{-6}$ \\
\hline
\end{tabular}
\endgroup
\end{table*}

Table~\ref{CoveragePDE} reports the reconstruction errors, variance-estimation errors, and empirical coverage. Debiasing reduces both mean reconstruction errors and mean absolute variance-estimation errors for all seven methods. Ordinary coverage is $100\%$ throughout, while debiased coverage ranges from $94.0\%$ to $97.1\%$, closer to the nominal $95\%$ level.

\FloatBarrier
\section{Conclusion and outlook}
\label{sec:Con}

This paper uses spectral analysis to establish a common statistical theory of debiasing, thresholding, and simultaneous inference for a class of linear regression estimators. The methods retain their own computational algorithms, including iterative updates and discretized dynamical equations. Generator and residual bounds provide a common analytical description of their errors. Debiasing uses two applications of the original linear algorithm with unchanged parameters; spectral analysis shows that debiasing replaces the residual by its square. Implementations based on matrix-vector products therefore retain their ability to avoid explicit spectral decomposition.

The statistical results establish maximum-norm consistency, recovery of the population threshold set, coefficient and noise-variance error bounds, and Gaussian approximations under explicit assumptions. The bootstrap theorem supplies simultaneous coverage guarantees under full column rank and deterministic sequences of regression parameters and thresholds. The separate best-worst-case analysis relates noise-free bias rates to perturbation error rates under additional regularity and gain conditions. Together, these results identify the spectral properties, signal restrictions, and tuning assumptions that connect the methods to their statistical guarantees.

The sparse and non-sparse regression experiments illustrate how debiasing and thresholding affect estimation accuracy under different signal structures and the specified tuning rules. The inverse source experiment further examines reconstruction accuracy, variance estimation, and empirical simultaneous coverage.

Two natural directions for future research are bootstrap inference with data-adaptive regularization, stopping, and threshold selection, and its extension to rank-deficient designs under suitable identifiability conditions.

\backmatter
\hypersetup{bookmarksdepth=2}
\bmhead{Supplementary information}
Online Resource~1 (\texttt{ESM\_1.pdf}) contains auxiliary probability
lemmas, the proofs of Theorems~\ref{Thm1}, \ref{them27},
\ref{them210}, \ref{thm:bootstrap-validity}, and~\ref{thmBW} and Lemma~\ref{BestWorst2},
method-specific spectral derivations and qualification bounds, detailed
numerical discretizations, additional simulation figures, and the
inverse-source formulation and reconstruction plots. The full proofs of
Theorems~\ref{thmP01}, \ref{thmP1}, \ref{them29}, and~\ref{them215}
are in Appendix~\ref{app:core-proofs} of this article.

\section*{Statements and Declarations}
\bmhead{Funding}
This work was supported by the Shenzhen Sci-Tech Fund (Grant No. RCJC20231211090030059), National Natural Science Foundation of China (Grant No. 12171036), National Key Research and Development Program of China (Grant No. 2022YFC3310300) and Beijing Natural Science Foundation (Grant No. Z210001).

\bmhead{Competing interests}
The authors declare no competing interests.

\bmhead{Data and code availability}
The code, numerical inputs, and saved results are available at
\url{https://github.com/Ao-King/HighDimLR}, version \href{https://github.com/Ao-King/HighDimLR/tree/fb2b57e497e0494a3646c6597dc4f6b6d40a297f}{\texttt{fb2b57e}}.

\bmhead{Author contributions}
Ying-Ao Wang developed the methodology, proved the theoretical results, implemented the algorithms, conducted the numerical experiments, analyzed the results, and drafted the manuscript. Yunyi Zhang and Ye Zhang conceived the study, supervised the research, and reviewed and revised the manuscript. Ye Zhang also provided research funding and equipment support.

\begin{appendices}
\allowdisplaybreaks[0]
\renewcommand{\theequation}{\thesection.\arabic{equation}}
\renewcommand{\theHequation}{\theHsection.\arabic{equation}}
\numberwithin{lemma}{section}
\renewcommand{\theHlemma}{\theHsection.\arabic{lemma}}
\section{Proofs of the core results}
\label{app:core-proofs}

We give the proofs of Theorems~\ref{thmP01}, \ref{thmP1},
\ref{them29}, and~\ref{them215}. The auxiliary moment and Gaussian approximation bounds are stated in
Supplementary Lemmas~\ref{supp-lemma1} and~\ref{supp-lemma:rectangle-clt}.

\subsection{Spectral decompositions and notation}
\label{app:spectral}

Use the thin singular value decomposition
\citep[Theorem 7.3.2]{johnson1985matrix}
\begin{equation}
\label{SVD}
\mathbf X_n=\mathbf U\varLambda\mathbf V^{\mathrm T}.
\end{equation}
Here $\varLambda=\operatorname{diag}\{\sqrt{\lambda_1},\ldots,
\sqrt{\lambda_s}\}$, $s=\operatorname{rank}(\mathbf X_n)\leq\min\{n,p\}$, and the
columns of $\mathbf U\in\mathbf R^{n\times s}$ and
$\mathbf V\in\mathbf R^{p\times s}$ are orthonormal. Thus
$\mathbf U^{\mathrm T}\mathbf U=\mathbf V^{\mathrm T}\mathbf V=\mathbf I_s$
and $\mathbf V\mathbf V^{\mathrm T}=\mathbf P_n$. With
$\pmb\beta_\perp=(\mathbf I_p-\mathbf P_n)\pmb\beta$, we have
\begin{equation}
\label{genest}
\begin{aligned}
\hat{\pmb\beta}_\alpha-\pmb\beta
&=\frac1n\mathbf Vg_\alpha(\varLambda^2/n)
\varLambda\mathbf U^{\mathrm T}\pmb e_n\\
&\quad-\mathbf Vr_\alpha(\varLambda^2/n)
\mathbf V^{\mathrm T}\pmb\beta-\pmb\beta_\perp,
\end{aligned}
\end{equation}
and
\begin{equation}
\label{debiasest1}
\begin{aligned}
\tilde{\pmb\beta}_\alpha-\pmb\beta
&=\frac1n\mathbf V[\mathbf I_s+r_\alpha(\varLambda^2/n)]\\
&\qquad\times g_\alpha(\varLambda^2/n)
\varLambda\mathbf U^{\mathrm T}\pmb e_n\\
&\quad-\mathbf Vr_\alpha^2(\varLambda^2/n)
\mathbf V^{\mathrm T}\pmb\beta-\pmb\beta_\perp.
\end{aligned}
\end{equation}

Write $\pmb\zeta=\mathbf V^{\mathrm T}\pmb\beta$,
$r_j=r_\alpha(\lambda_j/n)$, and $g_j=g_\alpha(\lambda_j/n)$, and set

\begin{equation*}
\begin{gathered}
B_i=\sum_{j=1}^s v_{ij}\zeta_jr_j,\qquad
\widetilde B_i=\sum_{j=1}^s v_{ij}\zeta_jr_j^2,\\
\ell_{il}=\frac1n\sum_{j=1}^s v_{ij}g_j\sqrt{\lambda_j}u_{lj},\\
\widetilde\ell_{il}=\frac1n\sum_{j=1}^s
v_{ij}(1+r_j)g_j\sqrt{\lambda_j}u_{lj},\\
Z_i=\sum_{l=1}^n\ell_{il}e_l,\qquad
\widetilde Z_i=\sum_{l=1}^n\widetilde\ell_{il}e_l.
\end{gathered}
\end{equation*}
The qualification condition and Assumptions 1--2 imply
\begin{equation*}
\rho_n:=\max_{1\leq j\leq s}|r_j|
\leq C_*\left(\frac{n\alpha}{\lambda_s}\right)^d
=O(n^{-d\delta})=o(1).
\end{equation*}
Consequently, for sufficiently large $n$,
\begin{equation*}
|g_j|\leq\frac{2n}{\lambda_j},
\qquad |(1+r_j)g_j|\leq\frac{4n}{\lambda_j}.
\end{equation*}
These absolute-value bounds justify the squared coefficient estimates
below. Write $q_n=\|\pmb\beta_\perp\|_\infty$; Assumption 5(b)(ii),
when imposed, gives $q_n=o((n\log n)^{-1/2})$.

\subsection{Consistency and threshold-set recovery}
\label{app:consistency}

\begin{proof}[Proof of Theorem~{\upshape\ref{thmP01}}]
\label{proof:thmP01}
By \eqref{genest},
$(\hat{\pmb\beta}_\alpha)_i-\beta_i
=Z_i-B_i-(\pmb\beta_\perp)_i$.
Since $\sum_jv_{ij}^2\leq1$ and
$\|\pmb\zeta\|_2\leq\|\pmb\beta\|_2$, Cauchy's inequality yields
\begin{equation}
\label{cauchy1}
\begin{aligned}
\max_i|B_i|&\leq\rho_n\|\pmb\beta\|_2\\
&\leq\frac{C_*(n\alpha)^d}{\lambda_s^d}\|\pmb\beta\|_2
=O(n^{\alpha_\beta-d\delta}).
\end{aligned}
\end{equation}
Orthonormality of the columns of $\mathbf U$ gives
\begin{equation}
\label{cauchy2}
\begin{aligned}
\max_i\sum_l\ell_{il}^2
&=\max_i\sum_j\frac{v_{ij}^2g_j^2\lambda_j}{n^2}\\
&\leq\max_i\sum_j\frac{4v_{ij}^2}{\lambda_j}
\leq\frac4{\lambda_s}.
\end{aligned}
\end{equation}
Supplementary Lemma~\ref{supp-lemma1}, with its constant $E$, implies
\begin{equation}
\label{eq29}
\mathbb P\left(\max_i|Z_i|>t\right)
\leq\frac{pE2^m}{\lambda_s^{m/2}t^m},\qquad t>0.
\end{equation}
Thus $\max_i|Z_i|=O_p(n^{\alpha_p/m-\eta})$. Also,
$\eta\leq1/2$ gives $q_n=o(n^{\alpha_p/m-\eta})$.
Since Assumption~2 gives $d\delta>\alpha_\beta+\eta$, \eqref{cauchy1}
is $o(n^{\alpha_p/m-\eta})$, and hence
\begin{align*}
\|\hat{\pmb\beta}_\alpha-\pmb\beta\|_\infty
&\leq\max_i|Z_i|+\max_i|B_i|+q_n\\
&=O_p(n^{\alpha_p/m-\eta}).
\end{align*}

For the debiased estimator, \eqref{debiasest1} gives
\begin{equation}
\label{IdentityBetaPoint}
(\tilde{\pmb\beta}_\alpha)_i-\beta_i
=\widetilde Z_i-\widetilde B_i-(\pmb\beta_\perp)_i.
\end{equation}
The corresponding bounds are
\begin{equation}
\label{cauchy3}
\max_i|\widetilde B_i|
\leq\rho_n^2\|\pmb\beta\|_2
=O(n^{\alpha_\beta-2d\delta})
\end{equation}
and
\begin{equation}
\label{cauchy4}
\max_i\sum_l\widetilde\ell_{il}^2
=\max_i\sum_j
\frac{v_{ij}^2(1+r_j)^2g_j^2\lambda_j}{n^2}
\leq\frac{16}{\lambda_s}.
\end{equation}
Supplementary Lemma~\ref{supp-lemma1} therefore gives
\begin{equation*}
\mathbb P\left(\max_i|\widetilde Z_i|>t\right)
\leq\frac{pE4^m}{\lambda_s^{m/2}t^m},\qquad t>0.
\end{equation*}
Assumption~2 also makes \eqref{cauchy3} negligible on the same scale.
Combining these bounds with \eqref{IdentityBetaPoint} and the bound on
$q_n$ proves the second rate.
\end{proof}

\begin{proof}[Proof of Theorem~{\upshape\ref{thmP1}}]
\label{proof:thmP1}
Let $h_b=\min\{c_b^{-1}-1,1-c_b\}>0$ and
$\varepsilon_b=h_b/3$. Assumption 2 gives
$\alpha_\beta-d\delta<-\eta$, while Assumption 4(b) gives
$v_b<\eta\leq1/2$. By \eqref{cauchy1}, \eqref{cauchy3}, and
Assumption 5(b)(ii),
\begin{equation*}
\max_i|B_i|+q_n=o(b_n),\qquad
\max_i|\widetilde B_i|+q_n=o(b_n).
\end{equation*}
For sufficiently large $n$, each left-hand side is at most
$\varepsilon_b b_n$. On $\max_i|Z_i|\leq\varepsilon_b b_n$,
\begin{equation*}
\max_i|(\hat{\pmb\beta}_\alpha)_i-\beta_i|
\leq2\varepsilon_b b_n<h_b b_n.
\end{equation*}
The separation in Assumption 4(b) then implies
$\widehat{\mathcal N}_{b_n}=\mathcal N_{b_n}$. By \eqref{eq29},
\begin{align*}
\mathbb P(\widehat{\mathcal N}_{b_n}\ne\mathcal N_{b_n})
&\leq\mathbb P(\max_i|Z_i|>\varepsilon_b b_n)\\
&\leq\frac{pE2^m}{\lambda_s^{m/2}(\varepsilon_b b_n)^m}
=O(n^{\alpha_p+mv_b-m\eta}).
\end{align*}
The same argument, using \eqref{cauchy4} and Supplementary
Lemma~\ref{supp-lemma1}, gives
\begin{equation*}
\mathbb P(\widetilde{\mathcal N}_{b_n}\ne\mathcal N_{b_n})
\leq\frac{pE4^m}{\lambda_s^{m/2}(\varepsilon_b b_n)^m}
=O(n^{\alpha_p+mv_b-m\eta}).
\end{equation*}
The exponent is negative by Assumption 4(b).
\end{proof}

\subsection{Gaussian approximation}
\label{app:gaussian}

Let $\ell_i$ and $\widetilde\ell_i$ denote the noise coefficient rows
defined above, and write $t_i^*=\ell_i/\tau_i^*$ and
$t_i=\widetilde\ell_i/\tau_i=a_i$. Then
$\tau_i^*=(\|\ell_i\|_2^2+1/n)^{1/2}$ and
$\tau_i=(\|\widetilde\ell_i\|_2^2+1/n)^{1/2}$.
For sufficiently large $n$, $\rho_n=o(1)$ gives
$|1-r_j|,|1-r_j^2|\geq1/2$. Under
\eqref{eq:gaussian-nondegeneracy}, there is a constant $c>0$ such that
\begin{equation}
\label{app:variance-bounds}
\begin{aligned}
\min\{\tau_i,\tau_i^*\}
&\geq\frac12\left(\sum_{j=1}^s
\frac{v_{ij}^2}{\lambda_j}\right)^{1/2},\\
c&\leq\min\{\|t_i\|_2^2,\|t_i^*\|_2^2\}\\
&\leq\max\{\|t_i\|_2^2,\|t_i^*\|_2^2\}\leq1.
\end{aligned}
\end{equation}
Indeed, $\|\widetilde\ell_i\|_2^2\geq c_a/n$ and
$\|\widetilde\ell_i\|_2\leq(1+\rho_n)\|\ell_i\|_2$.
If $s=p$, both unnormalized squared row norms are at least
$\tfrac14\sum_jv_{ij}^2/\lambda_j\geq1/(4C_\lambda n)$;
thus \eqref{eq:gaussian-nondegeneracy} holds automatically in that case.
Assumptions 1--2 give $d>1/2$ and
$\delta_1=d\delta-\eta-\alpha_\beta>0$. Weighted Cauchy--Schwarz
and the qualification condition yield
\begin{equation}
\label{app:normalized-bias}
\begin{aligned}
\max_i\frac{|B_i|}{\tau_i^*}
&\leq2\left(\sum_j\lambda_j\zeta_j^2r_j^2\right)^{1/2}\\
&\leq2C_*(n\alpha)^d\lambda_s^{1/2-d}\|\pmb\beta\|_2
=O(n^{-\delta_1}).
\end{aligned}
\end{equation}
For the debiased estimator the same argument gives
\begin{align}
&\max_i\frac{|\widetilde B_i|}{\tau_i}
\leq2\left(\sum_j\lambda_j\zeta_j^2r_j^4\right)^{1/2}\label{eq57}\\*
&\quad\leq2C_*^2(n\alpha)^{2d}
\lambda_s^{1/2-2d}\|\pmb\beta\|_2
=O(n^{-\delta_1-d\delta}).\notag
\end{align}
Moreover, Assumption 5(b)(ii) and $\tau_i,\tau_i^*\geq n^{-1/2}$ imply
\begin{equation}
\label{app:null-bias}
\max_i\frac{|(\pmb\beta_\perp)_i|}{\min\{\tau_i,\tau_i^*\}}
\leq\sqrt n\|\pmb\beta_\perp\|_\infty
=o((\log n)^{-1/2}).
\end{equation}

\begin{lemma}[Comparison of normalized noise rows]
\label{app:normalized-rows}
Under Assumptions 1--3, 4(a), and \eqref{eq:gaussian-nondegeneracy},
the rows defined above satisfy
\begin{equation}
\label{app:row-distance}
\max_i\|t_i^*-t_i\|_2
\leq\frac{2\rho_n}{1-\rho_n}=O(n^{-d\delta})
\end{equation}
for sufficiently large $n$. Let $J=J_n\subseteq\{1,\ldots,p\}$ be
any nonempty deterministic set, and let $e_l^*$ be independent
$N(0,\sigma^2)$ variables. Define
\begin{equation*}
A_J(z)=\max_{i\in J}\left|\sum_l t_{il}z_l\right|,\quad
A_J^*(z)=\max_{i\in J}\left|\sum_l t_{il}^*z_l\right|.
\end{equation*}
If
\begin{equation*}
\sup_{x\geq0}|\mathbb P(A_J(e)\leq x)
-\mathbb P(A_J(e^*)\leq x)|\longrightarrow0,
\end{equation*}
then the same convergence holds with $A_J^*$ in place of $A_J$.
\end{lemma}

\begin{proof}
Orthonormality of $\mathbf U$ gives
\begin{equation*}
\|\widetilde\ell_i-\ell_i\|_2^2
=\sum_j\frac{v_{ij}^2r_j^2g_j^2\lambda_j}{n^2}
\leq\rho_n^2\|\ell_i\|_2^2.
\end{equation*}
The reverse triangle inequality for the augmented vectors
$(\ell_i,n^{-1/2})$ and $(\widetilde\ell_i,n^{-1/2})$ gives
$|\tau_i-\tau_i^*|\leq\rho_n\|\ell_i\|_2\leq\rho_n\tau_i^*$.
In particular, $\tau_i\geq(1-\rho_n)\tau_i^*$. Consequently,
\begin{equation*}
\begin{aligned}
\|t_i-t_i^*\|_2
&\leq\frac{\|\widetilde\ell_i-\ell_i\|_2}{\tau_i}
+\frac{\|\ell_i\|_2|\tau_i-\tau_i^*|}{\tau_i\tau_i^*}\\
&\leq\frac{2\rho_n}{1-\rho_n},
\end{aligned}
\end{equation*}
which proves \eqref{app:row-distance}. Applying Supplementary
Lemma~\ref{supp-lemma1} to the difference rows, both for $e$ and
for $e^*$, yields
\begin{align*}
D_J(e)&:=\max_{i\in J}\left|\sum_l(t_{il}^*-t_{il})e_l\right|
=O_p(n^{\alpha_p/m-d\delta}),\\
D_J(e^*)&=O_p(n^{\alpha_p/m-d\delta}).
\end{align*}
Put $\gamma=d\delta-\alpha_p/m>0$ and $h_n=n^{-\gamma/2}$.
Then $\mathbb P(D_J(e)>h_n)+\mathbb P(D_J(e^*)>h_n)\to0$.
By \eqref{app:variance-bounds}, the Gaussian coordinates in
$A_J^*(e^*)$ have variances bounded above and away from zero.
The boundary bound in Supplementary equation~\eqref{supp-eq:rectangle-boundary}, which allows
singular covariance, gives, for their maximum's distribution function
$F_J^*$,
\begin{align*}
\omega_J(2h_n)&:=\sup_x\{F_J^*(x+2h_n)-F_J^*(x)\}\\
&=O(h_n\sqrt{\log n})=o(1).
\end{align*}
Let $\Delta_J$ denote the assumed approximation error for $A_J$.
The inequalities $|A_J^*(z)-A_J(z)|\leq D_J(z)$ give
\begin{align*}
&\sup_{x\geq0}|\mathbb P(A_J^*(e)\leq x)-F_J^*(x)|\\
&\quad\leq\Delta_J+\mathbb P(D_J(e)>h_n)\\
&\qquad+\mathbb P(D_J(e^*)>h_n)+\omega_J(2h_n).
\end{align*}
All terms tend to zero, proving the transfer statement.
\end{proof}

\begin{proof}[Proof of Theorem~{\upshape\ref{them29}}]
\label{proof:them29}
Use $J=\{1,\ldots,p\}$ in the preceding notation, and set
\begin{align*}
M_n&=\max_i\frac{|(\tilde{\pmb\beta}_\alpha)_i-\beta_i|}{\tau_i},\\
S_n&=A_J(e),\qquad S_n^G=A_J(e^*).
\end{align*}
By \eqref{debiasest1}, \eqref{eq57}, and \eqref{app:null-bias},
$|M_n-S_n|\leq\varepsilon_n$ for sufficiently large $n$, where
\begin{equation*}
\varepsilon_n=Cn^{-\delta_1}+\sqrt n\|\pmb\beta_\perp\|_\infty
=o((\log n)^{-1/2}).
\end{equation*}
Thus
\begin{equation}
\label{eq61}
\mathbb P(S_n\leq x-\varepsilon_n)
\leq\mathbb P(M_n\leq x)
\leq\mathbb P(S_n\leq x+\varepsilon_n).
\end{equation}
Since $t_{il}=a_{il}$, \eqref{eq:gaussian-nondegeneracy} gives
\begin{align}
\max_i\mathbb E\left(\sum_l t_{il}e_l^*\right)^2
&\leq\sigma^2,\label{eq63}\\
\min_i\mathbb E\left(\sum_l t_{il}e_l^*\right)^2
&\geq c_a\sigma^2>0.\label{eq64}
\end{align}
Write $a_{\max}=\max_{i,l}|a_{il}|$ and
$L_n=\sum_l\max_i|a_{il}|^3$. Under Assumption 6(A),
$L_n\leq n a_{\max}^3=o((\log n)^{-9/2})$.
Under Assumption 6(B), the same bound follows from
$L_n\leq a_{\max}\sum_{i,l}a_{il}^2\leq p a_{\max}$.
Since $p=O(n^{\alpha_p})$, Supplementary
Lemma~\ref{supp-lemma:rectangle-clt} yields
\begin{equation}
\label{eq68}
\Delta_n:=\sup_{x\geq0}|\mathbb P(S_n\leq x)
-\mathbb P(S_n^G\leq x)|\longrightarrow0.
\end{equation}
Orthonormality of $\mathbf U$ identifies the distribution of
$S_n^G$ with $H$. Supplementary equation~\eqref{supp-eq:rectangle-boundary} gives
\begin{align*}
\omega_H(\varepsilon_n)&:=\sup_x\{H(x+\varepsilon_n)-H(x)\}\\
&\leq C\varepsilon_n\sqrt{\log n}\longrightarrow0.
\end{align*}
The sandwich \eqref{eq61} bounds the required approximation error
by $\Delta_n+\omega_H(\varepsilon_n)\to0$.
\end{proof}

\begin{proof}[Proof of Theorem~{\upshape\ref{them215}}]
\label{proof:them215}
Take $J=N=\mathcal N_{b_n}$, which is nonempty for sufficiently large $n$.
Condition~\eqref{eq:gaussian-nondegeneracy} gives \eqref{app:variance-bounds}.
Restricting to these rows preserves the variance bounds and can only
decrease $L_n$ in the proof of Theorem~\ref{them29}. Thus
Supplementary Lemma~\ref{supp-lemma:rectangle-clt} gives
\begin{equation*}
\Delta_{n,N}:=\sup_{x\geq0}
|\mathbb P(A_N(e)\leq x)-H_N(x)|\longrightarrow0.
\end{equation*}
Here orthonormality of $\mathbf U$ identifies $H_N$ as the
distribution function of $A_N(e^*)$.

Set $M_{n,N}=\max_i{|\tilde\theta_i-\beta_i|}/\tau_i$ and let
$E_n=\{\widetilde{\mathcal N}_{b_n}=N\}$.
On $E_n$, the retained coordinates obey \eqref{debiasest1},
while every discarded coordinate has estimation error $-\beta_i$.
Since $\min\{\tau_i,\tau_i^*\}\geq n^{-1/2}$, Assumption 5(b)(i)
gives
\begin{align*}
q_{n,N}
&:=\max_{i\notin N}\max\left\{
\frac{|\beta_i|}{\tau_i},\frac{|\beta_i|}{\tau_i^*}\right\}\\
&\leq\sqrt n\max_{i\notin N}|\beta_i|
=o((\log n)^{-1/2}),
\end{align*}
with the maximum over an empty complement defined to be zero.
By \eqref{eq57} and \eqref{app:null-bias}, on $E_n$,
\begin{equation*}
\begin{aligned}
|M_{n,N}-A_N(e)|&\leq\varepsilon_n,\\
\varepsilon_n&:=Cn^{-\delta_1}+\sqrt n\|\pmb\beta_\perp\|_\infty+q_{n,N}\\
&=o((\log n)^{-1/2}).
\end{aligned}
\end{equation*}
Theorem~\ref{thmP1} gives
$\mathbb P(E_n^c)=O(n^{\alpha_p+mv_b-m\eta})=o(1)$.
The resulting sandwich implies
\begin{equation}
\label{eq710}
\begin{aligned}
&\sup_{x\geq0}|\mathbb P(M_{n,N}\leq x)-H_N(x)|\\
&\quad\leq\mathbb P(E_n^c)+\Delta_{n,N}
+\omega_{H_N}(\varepsilon_n),
\end{aligned}
\end{equation}
where $\omega_{H_N}(u)=\sup_x\{H_N(x+u)-H_N(x)\}$.
Supplementary equation~\eqref{supp-eq:rectangle-boundary} yields
\begin{equation}
\label{eq73}
\omega_{H_N}(\varepsilon_n)
\leq C\varepsilon_n\sqrt{\log n}\longrightarrow0.
\end{equation}

Equation~\eqref{eq710} therefore proves the first conclusion.

For the ordinary estimator, Lemma~\ref{app:normalized-rows},
applied on $N$, transfers the centered approximation to $A_N^*(e)$;
the corresponding Gaussian distribution is $H_N^*$.
On $E_n^*=\{\widehat{\mathcal N}_{b_n}=N\}$, the bounds
\eqref{app:normalized-bias}, \eqref{app:null-bias}, and the same tail bound give
\begin{equation*}
\left|\max_i\frac{|\hat\theta_i-\beta_i|}{\tau_i^*}
-A_N^*(e)\right|\leq\varepsilon_n.
\end{equation*}
Theorem~\ref{thmP1} gives $\mathbb P((E_n^*)^c)=o(1)$.
The variance bounds in \eqref{app:variance-bounds} and the same
boundary and sandwich argument apply to
$H_N^*$, proving the second conclusion.
\end{proof}
\end{appendices}

\makeatletter
\immediate\write\@auxout{\string\citation{SCBibliographySettings}}
\makeatother
\bibliography{references,journal-settings}
\end{document}